\documentclass[12pt]{article}
\usepackage[english]{babel}
\usepackage[a4paper,left=2.5cm,right=2.5cm,top=3cm,bottom=3cm]{geometry} 

\usepackage{mathtools}
\usepackage{bbm}
\usepackage[utf8]{inputenc}
\usepackage[T1]{fontenc}
\usepackage[english]{babel}
\usepackage{amstext,amssymb,amsopn,amsthm, amsfonts}
\usepackage{nicefrac, esint}
\usepackage{mathrsfs}
\usepackage{dsfont}
\usepackage{bbm}
\usepackage{thmtools}
\usepackage[colorinlistoftodos,prependcaption,color=white,bordercolor=blue]{todonotes}

\usepackage[shortlabels]{enumitem}
\usepackage[arrow, matrix, curve]{xy}

\usepackage{hyperref}
\hypersetup{
	colorlinks=true, 
	linktoc=all,     
	linkcolor=blue, 
	citecolor=black, 
	hypertexnames=false, %hack that seems to work idk 
}

\addto\extrasenglish{%
}

\usepackage{titlesec}
\usepackage{lipsum}
\titleformat{\section}
{\normalfont\scshape}{\thesection}{1em}{}
\titleformat{\subsection}
{\normalfont\normalfont\bfseries}{\thesubsection}{1em}{}
\titleformat{\subsubsection}[runin]
{\normalfont\normalsize\itshape}{\thesubsubsection}{1em}{}

\theoremstyle{plain}
\declaretheorem[title=Theorem, parent=section]{theorem}
\declaretheorem[title=Lemma,sibling=theorem]{lemma}
\declaretheorem[title=Proposition,sibling=theorem]{proposition}
\declaretheorem[title=Corollary,sibling=theorem]{corollary}

\theoremstyle{definition}
\declaretheorem[title=Definition,sibling=theorem]{definition}
\declaretheorem[title=Remark,sibling=theorem]{remark}
\declaretheorem[title=Remark, numbered=no]{remark*}
\declaretheorem[title=Example, sibling=theorem]{example}
\declaretheorem[title=Assumption, sibling=theorem]{assumption}

\numberwithin{equation}{section}

\newcommand{\borel}{\mathcal B}

\newcommand{\configMeas}{\nu^\Gamma}

\newcommand{\diamDom}{{\diam\domain}}

\newcommand{\domain}{\Omega}
\newcommand{\dPathMap}{\Phi}

\newcommand{\eps}{\varepsilon}

\newcommand{\fctempty}{\,\cdot\,}

\newcommand{\frontview}{\nu^\Gamma}
\newcommand{\fullview}{\mu^\Gamma}

\newcommand{\haus}{{\mathcal H}}
\newcommand{\id}{I}
\newcommand{\indic}{\mathbbm{1}}
\newcommand{\intd}{\mathrm{d}}
\newcommand{\inv}{^{-1}}

\newcommand{\leb}{\mathcal L}
\newcommand{\measConst}{\mu^{std}}

\newcommand{\N}{\mathbb{N}}

\newcommand{\pathP}{\mathbb{P}}
\newcommand{\pathMap}{\Psi}

\newcommand{\quadForm}{\mathcal E}

\newcommand{\R}{\mathbb{R}}

\newcommand{\Rd}{{\mathbb{R}^d}}
\newcommand{\rearview}{\rho^\Gamma}
\newcommand{\RotM}{A}
\newcommand{\Sd}{{\mathbb S^{d-1}}}

\newcommand{\thet}{\vartheta}
\newcommand{\upperInd}[1]{^{(#1)}}
\newcommand{\uxy}{\big(u(x)-u(y)\big)^2 \;|y-x|^{-d-2s} \;\intd y \intd x}
\newcommand{\vect}[1]{\mathbf{#1}}

\newcommand{\viewMap}{V}
\newcommand{\viewMapKT}{V^{(k)}_t}

\newcommand{\radAss}{{R_A}}
\newcommand{\lambdaAss}{{\lambda_A}}

\newcommand{\rhoN}{\rho^{(N)}}
\newcommand{\thetN}{\thet^{(N)}}
\newcommand{\phiN}{\phi^{(N)}}

\newcommand{\wurstN}{W^{(N)}}
\newcommand{\wurst}{W}

\newcommand{\constAss}{c_A}

\newcommand{\constCov}{\alpha}
\newcommand{\constDet}{c_D}

\newcommand{\constExp}{w}
\newcommand{\constL}{c_L}
\newcommand{\constI}{c_{\iota}}
\newcommand{\constM}{{M_w}}

\newcommand{\constW}{c_w}

\DeclareMathOperator{\adj}{adj}
\DeclareMathOperator{\conv}{conv}
\DeclareMathOperator{\diag}{diag}
\DeclareMathOperator{\diam}{diam}
\DeclareMathOperator{\dist}{dist}

\DeclareMathOperator*{\essup}{essup}
\DeclareMathOperator{\pv}{p.v.}

\DeclareMathOperator{\Tail}{Tail}

\makeatletter
\newcommand*{\defeq}{\mathrel{\rlap{%
			\raisebox{0.3ex}{$\m@th\cdot$}}%
		\raisebox{-0.3ex}{$\m@th\cdot$}}%
	=}
\makeatother

\numberwithin{equation}{section}

\begin{document}
	\title{Comparability for quadratic forms with space-dependent singular measures}
	
	\author{Kai-Uwe Bux and Solveig Hepp}
	\date{}
	\maketitle 
		
	\begin{abstract}
		\noindent
		This paper investigates quadratic forms $\quadForm^\Gamma$ generated by configurations $\Gamma$ and families of measures that are parametrized in the space variable and singular with respect to the Lebesgue measure. 
		A first basic result is the construction of corresponding integro-differential operators in divergence form. These operators are obtained by symmetrizing the family of singular measures.
		The main result is the comparability of the quadratic form $\quadForm^\Gamma$ with the quadratic form $\quadForm^s$ which is associated with the fractional Laplacian $(-\Delta)^s$.  
		Comparability of $\quadForm^\Gamma$ and $\quadForm^s$ implies Hölder regularity and a weak Harnack inequality for solutions to integro-differential equations involving these operators.
	\end{abstract}
	
	\renewcommand{\thefootnote}{\fnsymbol{footnote}} 
	\footnotetext{\emph{Date:} 8 November 2026} 	
	\renewcommand{\thefootnote}{\arabic{footnote}} 
	
\section{Introduction}\label{sec:intro}
	With a family of Borel measures $\mu(x, \intd y)$, parametrized in the space variable $x$ and possessing a singularity at $x = y$, one may formally associate an integro-differential operator
	\begin{align*}
		\mathcal L^\mu u(x) 
		\defeq \pv \; 2 \int_\Rd \big( u(x) - u(y)\big) \, \mu(x, \intd y). 
	\end{align*}
	If $\mu(x, \intd y)$ satisfies the symmetry condition
	\begin{align}\label{def:symmetry}
		\forall \; A, B \in \borel(\Rd) \; : \; 
		A \cap B = \emptyset
		\; \Rightarrow \;
		\int_A \int_B \mu(x, \intd y) \, \intd x 
		= \int_B \int_A \mu(x, \intd y)\,  \intd y, \tag{Sym}
	\end{align}
	we associate with $\mathcal L^\mu$ 
	a bilinear form 
	\begin{align*}%\label{def:quad-form-general-mu}
		\quadForm^\mu (u, v) 
		\defeq \int_\Rd \int_\Rd \big(u(x) - u(y)\big) \big(v(x) -v(y)\big) 
		\, \mu(x, \intd y) \,\intd x%\tag{$\quadForm^\mu$}
	\end{align*}
	via an integration by parts formula.
	In the study of integro-differential operators like $\mathcal L^\mu$, this bilinear form $\quadForm^\mu$ plays an important role. 
	On the one hand, in the area of stochastic processes, the existence of a strong Markov process with infinitesimal generator $\mathcal L^\mu$ follows from $\quadForm^\mu$ being a regular Dirichlet form on a suitable domain. 
	On the other hand, when investigating integro-differential equations like 
	\begin{align}\label{eq:intro-pde}
		\mathcal L^\mu u = 0 \text{ in } \domain \subset \Rd,\tag{IDE}
	\end{align}
	the concept of a weak solution may be defined using $\quadForm^\mu$. 
	
	The paradigmatic integro-differential operator is the fractional Laplacian 
	\begin{align}\label{def:fractional-laplacian}
		2 \, (- \Delta )^s u(x) 
		\defeq \text{ p.v. } \, c_{s, d} \, 2 \int_\Rd \big(u(x) - u(y)\big) |y-x|^{-d-2s} \, \intd y,
		\tag{FL}
	\end{align}
	with associated bilinear form 
	\begin{align*}
		\quadForm^s (u, v)
		= \int_\Rd \int_\Rd \big(u(x) -u(y)\big) \big(v(x) - v(y)\big) 
		\, |y-x|^{-d-2s} \, \intd y \intd x.
	\end{align*}
	Here, $2s$ is the operator's order of differentiability, with $s \in (0, 1)$. 
	
	Since much is known about $(-\Delta)^s$ and $\quadForm^s$, much can be learned by relating $\mathcal L^\mu$ and $\quadForm^\mu$ to $(-\Delta)^s$ and $\quadForm^s$. 
	An important tool in this context is the comparability of $\quadForm^\mu$ and $\quadForm^s$. We say that $\quadForm^\mu$ is \textit{comparable} to $\quadForm^s$ if there are constants $c_1, c_2 > 0$, such that, for all measurable functions $u: \Rd \to \R$,
	\begin{align*}
		c_1 \, \quadForm^s (u, u) 
		\leq \quadForm^\mu(u, u)
		\leq c_2 \, \quadForm^s (u, u).
	\end{align*}
	If $\quadForm^\mu$ satisfies the lower bound, it is called \textit{coercive}.
	The comparability of $\quadForm^\mu$ and $\quadForm^s$ has immediate consequences in the two research areas mentioned above. First, in the context of stochastic processes, the regularity of the Dirichlet form $\quadForm^s$ with the fractional Sobolev space $H^s(\Rd)$ as its domain almost implies the regularity of $\quadForm^\mu$ with the same domain.  
	Second, many regularity results for equations like \eqref{eq:intro-pde} assume coercivity of the quadratic forms associated with the operators $\mathcal L^\mu$. 
	
	When proving comparability, the choice of the family $\mu(x, \intd y)$ determines the difficulty of the task. 
	Absolutely continuous measures of the form $k(x, y) \, \intd y$, for example, have been investigated quite thoroughly. On the contrary, cases where $\mu(x, \intd y)$ is singular with respect to the Lebesgue measure are not yet very well understood. 
	A prominent example of such a family of singular measures is 
	\begin{align}\label{def:starrer-fall}
		\measConst(x, \intd y)
		%\defeq \sum_{i=1}^d \measConst_i (x, \intd y)
		\defeq \sum_{i=1}^d \Big[|y_i - x_i|^{-1-2s} \intd y_i \; \prod_{j \neq i} \delta_{\{x_j\}}(\intd y_j)\Big],\tag{$\measConst$}
	\end{align}
	where $\delta_{x_j}$ is the one-dimensional Dirac measure in the point $x_j$.
	This measure generates the operator $\sum_{i=1}^d (- \partial_i^2)^s$, the sum of one-dimensional fractional Laplacians. For its associated quadratic form, comparability can be proved by computing the respective operators' Fourier symbols. 
	
	The search for a comprehensive criterion on measures $\mu(x, \intd y)$ whose corresponding quadratic forms are coercive or even comparable to $\quadForm^s$ has gathered some attention. 
	One candidate is a nondegeneracy assumption, which requires $\mu(x, \intd y)$ to admit a constant $\lambda>0$, such that, for every point $x$, radius $r > 0$, and unit vector $\vect e \in \Sd$, there holds:
	\begin{align}\label{eq:conjecture-nondegeneracy}
		\int_{B_r(x)} \big((y-x) \cdot \vect e\big)^2_+ \; \mu(x, \intd y) 
		\geq \lambda r^{2-2s}\tag{Con}
	\end{align}
	Both \cite[(1.4)]{ChSi2020} and \cite[Condition (1.4)]{ImSi2020} conjecture that this nondegeneracy assumption implies comparability. 
	Recently, progress has been made on the conjecture, see comments below. 
	
	This paper addresses the question of comparability for bilinear forms governed by a class of measures $\mu(x, \intd y)$ where, for each $x$, the Borel measure $\mu(x, \fctempty)$ is singular with respect to the Lebesgue measure. 
	Unlike \eqref{def:starrer-fall}, however, the measures are allowed to vary in the space variable $x$.

\subsection{Main results}
	\label{sec:intro:main-results}
	We investigate families of singular measures that are constructed as follows. 
	\begin{definition}\label{def:configuration}
		A \textit{configuration} is a map $\Gamma = (\Gamma_1, \dots, \Gamma_d): \Rd \to (\Sd)^d$ consisting of \textit{subconfigurations} $\Gamma_k: \Rd \to \Sd, k \in \{1, \dots, d\}$, such that, for each $x \in \Rd$, $\Gamma(x)$ is a basis of $\Rd$, and each subconfiguration $\Gamma_k$ is $\borel(\Rd)-\borel(\Sd)$-measurable.
	\end{definition}
	Fixing $t \in \R$ and $k \in \{1, \dots, d\}$, define the \textit{view map} as
	\begin{align}\label{def:view-map}
		\viewMapKT: \Rd &\to \Rd \tag{VM}\\
		x&\mapsto x + t \Gamma_k(x),\notag
	\end{align}
	which maps a point $x$ to the point on the line $x + \langle \Gamma_k(x)\rangle$ that is $t$-far away from it. 
	%The \textit{nonsymmetric class of singular measures} $\frontview(x, \intd y)$ is given by:
	\begin{definition}\label{def:frontview-concrete}
		Let $T \in (0, \infty)$. Let $\Gamma: \Rd \to (\Sd)^d$ be a configuration such that the differential $D \viewMapKT$ exists almost everywhere for any $t \in (-T, T)$.
		Let $\iota: \leb(\R) \to [0, \infty]$ be a one-dimensional Lebesgue measure such that 
		\begin{align}\label{def:iota-comparable}
			\iota (\intd t) = I(t) \, \intd t  
			\text{ with }
			\exists \; \constI \geq 1 \; \text{ : } \; 
			\constI\inv |t|^{-1-2s} 
			\leq I(t) 
			\leq \constI |t|^{-1-2s}
			\tag{1D}
		\end{align}
		for $s \in (0, 1)$ and $t \in (-T, T)$. 
		The \textit{family of singular measures generated by} $\Gamma$ \textit{and} $\iota$ \textit{and cut off at} $T$ is given by
		\begin{align*}
			\nu^{\Gamma}(x, \fctempty) : \borel(\Rd) &\to [0, \infty]\\
			A &\mapsto \sum_{k=1}^d \int_{-T}^T
			\indic_A \big( x + t \Gamma_k (x)\big) 
			|\det D \viewMapKT(x)| 
			\, \iota(\intd t),
		\end{align*}
		where $x \in \Rd$.
		In this case, the configuration $\Gamma$ is called the \textit{generating configuration} of the family of measures $\nu^{\Gamma} (x, \intd y)$, and $T$ its \textit{cutoff point}.
	\end{definition}
	
	We may think about these measures as follows. The configuration $\Gamma$ attaches a basis $\Gamma(x)$ to every point $x \in \Rd$. Then $\frontview (x, \fctempty)$ measures the intersections of the lines $x + \langle\Gamma_k(x)\rangle$ extending $T$-far from $x$ with the Borel set $A$ using the weighted Lebesgue measure $|\det D \viewMapKT(x)| \, \iota(\intd t)$. 
	We may imagine that every point "sees" along each of the $d$ vectors belonging to its basis $\Gamma(x)$, and that $\frontview(x, \intd y)$ measures each intersection of these axes of sight with the set $A$.

	For the family of measures $\frontview(x, \intd y)$ given in \autoref{def:frontview-concrete} and an open set $\domain \subset \Rd$, we define a symmetric bilinear form via
	\begin{align}\label{def:quad-form-gamma}
		\quadForm^\Gamma_\domain (u, v) 
		\defeq \int_\domain \int_\domain \big(u(x) - u(y)\big) \big(v(x) - v(y)\big)
		\, \frontview(x, \intd y) \, \intd x.
		\tag{$\quadForm^\Gamma$}
	\end{align}
	Even though $\quadForm^\Gamma_\domain$ is always symmetric in $u$ and $v$, the \textit{family of measures} $\frontview(x, \intd y)$ is, in general, {not} symmetric in the sense of \eqref{def:symmetry}, see \autoref{meas:ex:nonsymmetric-measure}. It therefore does not give rise to an integro-differential equation \eqref{eq:intro-pde} in divergence form.\footnote{We use the naming conventions of divergence and non-divergence form nonlocal operators established in \cite{FRROKorea}, see also \cite{FRROBook}.}
	For $s \in (0, 1)$, let 
	\begin{align*}%\label{def:quad-form-s}
		\quadForm^s_\domain (u, u) 
		\defeq \int_\domain \int_\domain \big( u(x) - u(y) \big)^2 
		\,|y-x|^{-d-2s} \, \intd y \intd x. %\tag{$\quadForm^s$}
	\end{align*}
	This is the quadratic form associated with the fractional Laplacian, restricted to the Cartesian square of $\domain$.
	
	In order to prove our main result -- the comparability of $\quadForm^\Gamma_\domain$ and $\quadForm^s_\domain$ -- we need to impose some restrictions on the behavior of the configuration $\Gamma$. Most importantly, we will assume:
	\begin{align}
		&\tag{Lip}\label{ass:lip}
		\text{Each subconfiguration }\Gamma_k \text{ is Lipschitz continuous with constant }\constL > 0.
	\end{align}
	This lets us control the configuration's behavior on small scales. In order to control its behavior on large scales, we will assume:
	\begin{align}\label{ass:w}
		\text{There are } \constM \in \N, \constExp \in [0, 2s), \constW > 0, \text{ such that, for every } k \in \{1, \dots, d\} \notag\\
		\text{and a.e. } y \in \Rd, \text{ there holds } \haus^0\big(\{x \in \Rd \mid \viewMapKT(x) = y\}\big) 
		\leq \max\{\constM, \constW |t|^\constExp\},\tag{W}
	\end{align}
	where $\haus^0$ is the zero-dimensional Hausdorff (or counting) measure. 
	This condition says that the cardinality of the preimage set $\{ x \in \Rd \mid \viewMapKT(x) = y\}$ scales with $|t|$ if $|t|$ is large. 
	A third assumption concerns the interaction between the subconfigurations $\Gamma_1, \dots, \Gamma_d$.
	If $\RotM_x$ is the matrix containing the vectors $\Gamma_1(x), \dots, \Gamma_d(x)$ as its columns, assume: 
	\begin{align}
		&\tag{Det}\label{ass:det}
		\text{There is a constant }\constDet \in (0, 1) \text{ such that } \constDet < |\det \RotM_x |.
	\end{align}
	Put differently, assumption \eqref{ass:det} says that there is a minimal (positive) angle between any two vectors of each basis $\Gamma(x) = \big(\Gamma_1(x) , \dots, \Gamma_d(x)\big)$. These three assumptions are discussed in more detail in \autoref{sec:measures:configs}.
	
	We prove:  
	\begin{theorem}[Comparability for Lipschitz configurations]\label{thm:comparability-lip}
		Let $s,\constDet  \in (0, 1)$, $\constI > 1$, $\constL, \constW, T_0 > 0$, $\constExp \in [0, 2s)$, and $\constM \in \N$. Let $\domain \subset \Rd$ be an open, bounded, and connected domain with Lip\-schitz boundary. 
		There are constants 
		\begin{align*}
			c_1 = c_1\big(d, s, \domain, \constI, \constL, \constDet, T_0\big) > 0
			\text{ and }
			c_2 = c_2 \big(d, s, \domain, \constI, \constL, \constM, \constExp, \constW\big) > 0, 
		\end{align*}
		such that, for any measurable function $u: \domain \to \R$, there holds 
		\begin{align}\label{eq:comparability}
			c_1\, \quadForm^s_\domain (u, u) 
			\leq \, \quadForm^{\Gamma}_\domain(u, u)
			\leq c_2 \,\quadForm^s_\domain (u, u),
			\tag{V}
		\end{align}
		where $\quadForm^\Gamma_\domain$ is a quadratic form with cutoff point $T \in (T_0, \infty]$ generated by 
		a configuration $\Gamma: \Rd \to (\Sd)^d$ 
		that satisfies
		\eqref{ass:lip} with Lipschitz constant $\constL$,
		\eqref{ass:det} with constant $\constDet$, 
		and \eqref{ass:w} with constants $\constM, \constExp$, and $\constW$, 
		and by a one-dimensional Lebesgue measure $\iota(\intd t)$ that satisfies \eqref{def:iota-comparable} with $\constI$ for $s$.
	\end{theorem}
	The dependencies of the constants show that assumption \eqref{ass:w} is only needed for the upper bound of \eqref{eq:comparability}, while \eqref{ass:det} is required solely for the lower bound. 
	Let us give an intuition on why this is so by sketching the method of proof for the lower bound in \eqref{eq:comparability}. 
	Intuitively, the difficulty of this inequality lies in the fact that $\quadForm^s_\domain$ integrates over all tuples $(x, y) \in \domain \times \domain$, while $\quadForm^\Gamma_\domain$ only takes into account tuples of the form $\big( x, x + t \Gamma_k(x)\big)$. 
	The remedy for this is to relate each tuple $(x, y) \in \domain \times \domain$ uniquely to a chain of tuples 
	\begin{align}\label{intro:def:chain-of-tuples}
		\big(x = p_0(x, y), p_1(x, y)\big), \;
		\big(p_1(x, y), p_2(x, y)\big), \;
		\dots, \;
		\big(p_{N-1}(x, y), p_N(x, y) = y\big)
	\end{align}
	which appear in the double integral $\quadForm^\Gamma_\domain$. If the distances $|p_n(x, y) - p_{n-1}(x, y)|$ are comparable to $|x-y|$ as well, one may hope to achieve an estimate like 
	\begin{align*}
		\big(u(x) - &u(y)\big)^2 \;|y-x|^{-d-2s}\\
		&\leq c
		\sum_{n=1}^N \big( u(p_{n-1}(x, y)) - u(p_n(x, y))\big)^2 \; \big|p_{n}(x, y) - p_{n-1}(x, y)\big|^{-1-2s},
	\end{align*}
	from which the lower bound of \eqref{eq:comparability} would follow.
	
	A tuple $(x, y)$ appears in the double integral $\quadForm^\Gamma_\domain(u, u)$ if there are $t \in \R$ and $k \in \{1, \dots, d\}$, such that $x + t \Gamma_k(x) = y$ or $y + t \Gamma_k (y) = x$. If this is the case, we will say that $x$ and $y$ are \textit{directly connected}. If there is a chain of points $x =: p_0(x, y), p_1(x, y), \dots, p_N(x, y) \defeq y$, such that two consecutive points are directly connected, we say that $x$ and $y$ are \textit{connected under the configuration} $\Gamma$. 
	Clearly, only if all points $x, y \in \domain$ are connected under $\Gamma$ may one hope to find a chain of tuples \eqref{intro:def:chain-of-tuples}. 
	The identification of classes of configurations that allow such a connectedness result is another topic of this paper that is also of independent interest. 
	
	%An assumption like \eqref{ass:det} ensures that the axes of sight $\Gamma_1(x), \dots, \Gamma_d(x)$ in each point not only span the whole space but also cannot converge to a lower-dimensional subspace of $\Rd$. This prevents the potential occurrence of dead ends, where groups of points cannot reach certain areas of the space, which would make a connectedness result impossible. 
	%Assumption \eqref{ass:w}, on the other hand, makes sure that $\quadForm^\Gamma_\domain$ cannot get too large by disallowing sets which are seen by infinitely many other sets, which is relevant for the upper bound in \eqref{eq:comparability}. 
	
	To the authors' best knowledge, \autoref{thm:comparability-lip} is the first comparability result that deals with non-translation invariant singular measures. 
	
	\medskip
	
	Equations in divergence form like \eqref{eq:intro-pde} do not in general admit pointwise solutions, since $\mathcal L^\mu u(x)$ does not exist pointwise -- even for very smooth functions. Therefore, when investigating equations like \eqref{eq:intro-pde}, one establishes a concept of weak solution using the bilinear form $\quadForm^\mu$, see e.g. \cite[Section 2.2.4]{FRROBook} or \cite[Definition 1.2]{Chen2024}. This bilinear form is formally associated with the operator $\mathcal L^\mu$  via an integration by parts formula that presupposes the symmetry of $\mu(x, \intd y)$. 
	
	As has been pointed out, the family $\frontview(x, \intd y)$ is not symmetric. 
	We thus construct a family of measures $\fullview(x, \intd y)$ satisfying the symmetry condition \eqref{def:symmetry} from $\frontview(x, \intd y)$. 
	%This then allows the definition of weak solutions to integro-differential equations like \eqref{eq:intro-pde} for operators in divergence form. 
	In \autoref{thm:properties-for-regularity-theory}, we establish properties of this symmetric family of measures $\fullview(x, \intd y)$ which, together with \autoref{thm:comparability-lip}, lead to Hölder regularity and a weak Harnack inequality for weak solutions to equations in divergence form governed by $\fullview(x, \intd y)$.

	We show: 
	\begin{theorem}\label{thm:properties-for-regularity-theory}
		Let $\constL> 0$. There exists some $T_0 = T_0(d, \constL) > 0$, such that
		every configuration $\Gamma: \Rd \to (\Sd)^d$ that satisfies \eqref{ass:lip} with Lipschitz constant $\constL$, 
		every one-dimensional Lebesgue measure  $\iota(\intd t)$ that satisfies \eqref{def:iota-comparable} for some $s \in (0, 1)$,
		and every cutoff point $T \in (0, T_0)$ generate a family of Borel measures 
		\begin{align} \label{def:fullview-concrete}
			\fullview (x, A) 
			\defeq 
			\sum_{k=1}^d \frac 12
			\Big[
			&\int_{-T}^T \indic_A \big(x + t \Gamma_k(x) \big) |\det D \viewMapKT(x)| \, \iota (\intd t) \tag{$\fullview$}\\
			+ &\int_{-T}^T \haus^0 \big( \{ y \in A \mid \viewMapKT(y) = x\}\big)
			\, \iota(\intd t)
			\Big],
			\notag
		\end{align} 
		where $x \in \Rd$ and $A \in \borel(\Rd)$,
		which is well-defined in the sense that 
		\begin{align}\label{def:well-definedness-measures} 
			\forall \; B \in \borel(\Rd) \; : \;
			x \mapsto \fullview(x, B)
			\text{ is } \borel(\Rd)-\borel([0, \infty])\text{-measurable.}\tag{WD}
		\end{align}
		Moreover, there holds: 
		\begin{enumerate}[(a)]
			\item The family of measures $\fullview(x, \intd y)$ is symmetric in the sense of \eqref{def:symmetry}. 
			\item If $\Gamma$ additionally satisfies \eqref{ass:w}, then $\fullview(x, \intd y)$ satisfies a L\'evy-like condition  
			\begin{align}\label{def:levy}
				\essup_{x \in \Rd} \int_\Rd \min\{1,|y-x|^2\} \, \fullview(x, \intd y) < \infty. \tag{L}
			\end{align}
			and, for $s$, the \textit{decay condition} 
			\begin{align}\label{def:decay}
				\exists \; c > 0 \; 
				\text{ : } \, \forall\; r > 0 \; \: \; 
				\essup_{x_0 \in \Rd} \mu \big(x_0, B_r(x_0)^c\big)
				\leq c r^{-2s}.\tag{Dec}
			\end{align} 
			\item If $\Gamma$ additionally satisfies \eqref{ass:w}, then there holds, for every open and bounded set $\domain \subset \Rd$ and any measurable function $u: \domain \to \R$,
			\begin{align*}
				\quadForm_\domain^\Gamma(u, u) = \quadForm_\domain^{\fullview}(u, u) 
				= \int_\domain \int_\domain \big(u(x) - u(y)\big)^2 \fullview(x, \intd y) \intd x.
			\end{align*}
		\end{enumerate}
	\end{theorem}
	
	We show in \autoref{sec:measures} that not all configurations $\Gamma$ and one-dimensional measures $\iota(\intd t)$ allow the definition of a family of measures $\fullview(x, \intd y)$ with the properties stated in the theorem. 
	In \autoref{meas:ex:non-symmetrizable-family-of-measures}, we construct a configuration that does not give rise to a symmetric family of measures $\fullview(x, \intd y)$. 
	Lipschitz continuity, however, always guarantees the existence of such a family of measures. The weight $|\det D \viewMapKT(x)|$ is needed precisely in this construction. 	
	The decay behavior \eqref{def:decay} is implied by assumption \eqref{ass:w}, see also \autoref{measures:ex:decay-counterexample}.
		
	\medskip 
	
	In recent years, there has been significant interest in integro-differential equations in divergence form like \eqref{eq:intro-pde}. 
	%	\begin{align}\label{appl:eq:pde-div-form}
		%		\mathcal L^\mu u(x) 
		%		\defeq \pv \int_\Rd \big(u(x) - u(y)\big) \, \mu(x, \intd y)
		%		= 0 \text{ for } x \in \domain \subset \Rd.\tag{DFE}
		%	\end{align}
	%These equations do not, in general, admit strong solutions, since $\mathcal L^\mu u(x)$ does not exist pointwise, even for very smooth functions. 
	One line of question is which assumptions on the family of measures $\mu(x, \intd y)$ still allow regularity or a weak Harnack inequality for weak solutions. 
	A very general result in this direction is \cite{Chen2024}, where  integro-differential operators defined with families of measures $\mu(x, \intd y)$ that must be neither absolutely continuous with respect to the Lebesgue measure nor translation invariant are studied.  
	Under the mere assumption of well-definedness \eqref{def:well-definedness-measures}, symmetry \eqref{def:symmetry}, appropriate decay \eqref{def:decay}, and coercivity of the quadratic form, the author shows Hölder regularity and a weak Harnack inequality for weak solutions to \eqref{eq:intro-pde}, see \cite[Assumptions 1.1 - 1.4 and Remark 1.1, as well as Theorems 1.1 and 1.2]{Chen2024}.
	
	Using these theorems together with \autoref{thm:properties-for-regularity-theory} and \autoref{thm:comparability-lip}, we infer: 
	\begin{theorem}[Hölder estimates and weak Harnack inequality for non-translation invariant singular measures]
		\label{thm:regularity}
		Let $\constL > 0$ and $s \in (0, 1)$. There is a cutoff point $T_0 = T_0(x, \constL) > 0$, 
		such that 
		any configuration $\Gamma: \Rd \to (\Sd)^d$ satisfying \eqref{ass:lip} with $\constL$, \eqref{ass:w}, and \eqref{ass:det}, 
		any one-dimensional Lebesgue measure $\iota(\intd t)$ satisfying \eqref{def:iota-comparable} for $s$,
		and any cutoff point $T \in (0, T_0)$, 
		define a symmetric family of singular measures parametrized in $x \in \Rd$ that is given by \eqref{def:fullview-concrete}.
		
		For the corresponding integro-differential equation \eqref{eq:intro-pde}, where $\mu = \fullview$, there holds: 
		\begin{enumerate}[(a)]
			\item Any locally bounded weak solution $u$ of \eqref{eq:intro-pde} is locally Hölder continuous in $\domain$.
			\item For any weak supersolution of \eqref{eq:intro-pde} that additionally satisfies $u \geq 0$ on some ball $B_{R}(x_0)$, a weak Harnack inequality holds, i.e. for any $r \in \big( 0, \frac{R}{42}\big)$
			\begin{align*}
				\inf_{B_{2r}(x_0)} u 
				\geq c \Big(\fint_{B_r(x_0)} u^\eps(x) \, \intd x\Big)^{\frac{1}{\eps}}
				- \Big(\frac{r}{R}\Big)^{2s} \Tail(u_-; x_0, \tfrac{7}{8} R, R),
			\end{align*}
			where $c, \eps \in (0, 1)$ do not depend on $u$. 
		\end{enumerate}	
	\end{theorem} 
	
	\begin{proof}
		The existence of the symmetric family of singular measures $\fullview(x, \intd y)$ given by \eqref{def:fullview-concrete} follows from \autoref{thm:properties-for-regularity-theory}, which also provides properties \eqref{def:symmetry}, \eqref{def:well-definedness-measures}, and \eqref{def:decay}. 		
		Moreover, by \autoref{thm:properties-for-regularity-theory} (c), $\quadForm^{\fullview}_\domain(u, u)$ and $\quadForm^\Gamma_\domain(u, u)$ coincide on measurable functions. 
		Hence, \autoref{thm:comparability-lip} gives the coercivity of the quadratic form $\quadForm^{\mu^\Gamma}_\domain$, which means that all assumptions on the family of measures $\mu(x, \intd y)$ required for \cite[Theorems 1.1 and 1.2]{Chen2024} are met. This implies Hölder continuity of weak solutions of \eqref{eq:intro-pde} and the weak Harnack inequality for weak supersolutions. 
	\end{proof}
	
	While the scope of \cite{Chen2024} is very general, so far, it was not obvious that the class of measures satisfying the four necessary assumptions contains any interesting examples that have not already been covered under more restrictive assumptions such as absolute continuity or translation invariance. The symmetric family of singular measures $\fullview(x, \intd y)$ provides such an example. 
	This allows regularity results for a novel class of operators in divergence form.
	
\subsection{Related literature}\label{sec:intro:related-literature}
	As mentioned above, conjecture \eqref{eq:conjecture-nondegeneracy} has attracted attention for several years. 
	Substantial progress has recently been made  by Behn, Colombo, Fern\'andez-Real, Grube, Prasad, Ros-Oton and will be published shortly. 
	An AI-based counterexample can be found in \cite{Biswas2026}. 
	
	Let us now provide an overview of the literature on comparability results. It is clear that comparability holds for measures $k(x, y) \, \intd y$ with $k(x, y) \asymp |y-x|^{-d-2s}$. 
	One way to generalize this observation is to restrict the area where this kind of pointwise comparability of the density function $k$ holds. 
	An early result in this direction is \cite{Dyda2006}. The restriction is effected by a function $h: \Rd \to [0, \infty]$ which satisfies $h(rz) = h(z)$ for all $z \in \Rd$, $r > 0$, and is bounded away from zero on a non-empty open set. The author shows
	\begin{align*}
		\int_\domain \int_\domain &\big( u(x) - u(y)\big)^2 |y-x|^{-d-2s} \, \intd y \intd x\\
		&\leq c \int_\domain \int_\domain 
		\big( u(x) - u(y)\big)^2 h(x-y)
		|y-x|^{-d-2s} \, \intd y \intd x.
	\end{align*}
	Here, $\domain$ is a bounded, connected Lipschitz or other similarly well-behaved domain. 
	An illustrative example is $h(x-y) = \indic_V(x-y)$, where $V$ is some double cone in $\Rd$.
	Since $h$ essentially causes a restriction of the domain of integration, the strategy of the proof is to find, for each pair $x, y$, a sequence of points $x = p_0, p_1, \dots, p_n = y$, such that $h(p_{n-1} - p_{n})$ is bounded from below. 
	Note that, in \cite{Dyda2006}, $h(x-y)$ can even be replaced by the function $\indic_{B_{\frac 12 \dist(x, \partial \domain)}(x)}(y)$. 
	Similar results may be found in \cite[Theorems 1.5, 1.6]{PratsSaks2017}. 
	
	The result of \cite{Dyda2006} can be generalized by allowing the area of pointwise comparability of $k(x, y)$ to vary in $x$. This idea is implemented, for example, in \cite{BKS2019} and \cite{ImSi2020}.
	In \cite{BKS2019}, kernels $k(x, y)$ comparable to $\big(\indic_{V[x]}(y) + \indic_{V[y]}(x)\big)|y-x|^{-d-2s}$ are considered, where $V[x]$ is a double cone assigned to each point $x$ in the space by a measurable function. The only requirement on these cones is that their apex angles  be bounded uniformly from below. In other words, the cones are not allowed to degenerate to a line. The authors show a comparability result like \eqref{eq:comparability} with a  chaining technique analogous to \cite{Dyda2006}. 
	%For each tuple $(x, y) \in \Rd \times \Rd$, a point $z$ is constructed, which lies in the intersection $V[x] \cap V[y]$. Then the difference $u(x) - u(y)$ is estimated via this middle point. 
	In a similar way, \cite[Lemma A.6]{ImSi2020} discusses kernels that are comparable on bundles of rays which may vary in $x$. While these bundles do not need to be double cones as in \cite{BKS2019}, they still always have positive Lebesgue measure in $\Rd$. The exact properties of these kernels come from the Boltzmann equation. 
	The technique of inserting chains of points to make up for gaps in the area of integration is explicitly baked into the assumptions on the kernel for whose quadratic form a comparability estimate is shown in \cite{HuKa2007}.  
%	The authors define for each $n \in \N$ a partition of $\Rd$ in the form of $\mathcal Q_n \defeq \bigcup_{x \in \frac{1}{n}\Z^d} \{\bigtimes_{i=1}^d [x_i, x_i + \frac{1}{n})\}$. They consider kernels $k(x, y)$, for which there is a number $K \in \N$ and a constant $\lambda > 0$, such that, for any $\eps > 0$ and  $n \in \N$, two cubes $Q, Q' \in \mathcal Q_n$ can be connected by a chain $Q = Q_0, Q_1, \dots, Q_K = Q'$ of other $\mathcal Q_n$-cubes, such that 
%	\begin{align*}
%		\int_{Q}\int_{Q'} \indic_{\{|y-x|>\eps\}}(x, y) |y-x|^{-d-2s} \, \intd y \intd x
%		\leq \lambda \int_{Q_{k-1}} \int_{Q_k} \indic_{\{|y-x|>\eps\}}(x, y) k(x, y) \, \intd y \intd x,
%	\end{align*}
%	for any $k \in \{1, \dots, K\}$.
%	They also need an assumption preventing the same cubes from being used too often in these chains. By inserting points lying in these intermediate cubes $Q_k$ in a similar fashion to the approach in \cite{Dyda2006}, the authors deduce a comparability result for quadratic forms with kernels of this kind. 
	
	The most general result working with restrictions of the domain of pointwise comparability is \cite{ChSi2020}. 
	The authors discuss families of measures $k(x, y) \, \intd y$, for which there exist constants $c_1 \in (0, 1), c_2 > 0$, such that, for any ball $B \in \Rd$ and $x \in B$, there holds $|\{ z \in B \mid k(x, z) \geq c_2 |x-z|^{-d-2s}\}| \geq c_1 |B|$ (see their Assumption 1.1). Intuitively, this captures the requirement that the set where the kernel $k(x, y)$ is comparable to $|x-y|^{-d-2s}$ has some density in all directions. Measures defined via these kernels allow a proof of the lower bound in \eqref{eq:comparability}. Interestingly, the result is not proved via a chaining technique. Instead, the authors construct a sequence of comparable kernels whose comparability sets grow, similar to the growing ink spots lemma due to Krylov and Safonov. 
	The result in \cite{ChSi2020} contains the results of \cite{Dyda2006}, \cite{HuKa2007}, \cite{BKS2019}, and \cite{ImSi2020}. All measures considered in these results are absolutely continuous with respect to the Lebesgue measure. 
	
	In \cite{DyKa20}, comparability for some singular measures is discussed. In particular, \eqref{eq:comparability} is shown for the family of measures \eqref{def:starrer-fall}, which corresponds to the family $\frontview(x, \intd y)$ generated by the constant configuration $\Gamma \equiv (\vect e_1, \dots, \vect e_d)$ in our setting. 
	While the authors are able to treat quite general measures, including singular ones, they must make the rather strong assumption of near-translation invariance: The families of measures $\mu(x, \intd y)$ must satisfy, for every $x \in \Rd$ and $f(x, \fctempty): \Rd \to [0, \infty]$ measurable,
	\begin{align*}
		\int_\Rd f(x, x+z) \nu_\ast(\intd z)
		\leq \int_\Rd f(x, y)\mu(x, \intd y)
		\leq \int_\Rd f(x, x+z) \nu^\ast(\intd z),
	\end{align*}
	where $\nu^\ast$ and $\nu_\ast$ are Borel measures (cf. \cite[Assumption (T)]{DyKa20}). They then show comparability for the quadratic form belonging to the constant case by a chaining argument. The families of measures $\frontview(x, \intd y)$ we consider do not require such an assumption, which makes \autoref{thm:comparability-lip} a more general result than the one of \cite{DyKa20}.

	\medskip
	
	Next, we will review related results on regularity theory. 
	
	Hölder estimates for elliptic equations in divergence form were first established by \cite{Kasmann2009} for stable-like operators. A recent result in the same setting is \cite{Fall2020}. With an additional regularity assumption on the kernel, \cite{Fall2020} is able to recover the intuition of solutions gaining $2s$ derivatives. 
	Stable-like operators in divergence form have extensively been treated in the parabolic setting as well. Hölder estimates and a weak Harnack inequality are proved in \cite{Komatsu1988}, \cite{Komatsu1995}, \cite{ChKu2003}, \cite{ChKu2008}, \cite{CaEtAl2011}, \cite{DiCaEtAl2014}, and \cite{DiCaEtAl2016}.
	In \cite{BaBaCheKa2009} and \cite{BaKaKu2010}, the stable-like setting is slightly modified to allow different jump intensities, i.e. the kernel $k(x, y)$ is pointwise comparable from above and from below to $|y-x|^{-d-\beta_i}$ with $i=1, 2$. In \cite{BaBaCheKa2009}, a Harnack inequality is established as well as an example of a solution that is not continuous. Under an additional weak smoothness assumption on the kernel, Hölder regularity is recovered in \cite{BaKaKu2010}. 
	
	More general kernels are often treated under the assumption of comparability. In the parabolic setting, for example, \cite{FeKa2013} and \cite{KaSchwa2014} establish Hölder regularity and a weak Harnack inequality. \cite{KaSchwa2014} can even treat singular kernels. 
	In the elliptic setting, one of the most general results is \cite{DyKa20}, where only the existence of cutoff functions is required in addition to comparability. In this work, the authors give explicit examples for the singular kernels they are able to treat, which are linear transformations of the constant family of singular measures \eqref{def:starrer-fall}. 
	Similarly, \cite{ChaKa2020} treats singular kernels like \eqref{def:starrer-fall} whose orders of differentiability are allowed to vary in each coordinate. It does not assume any kind of smoothness in the $x$-variable, but requires a comparability result with the quadratic form associated with \eqref{def:starrer-fall}. 
	\cite{FRROKorea} establishes a priori Hölder estimates for kernels which may be singular as well. However, their kernels must be almost even -- a symmetry assumption that is more natural to non-divergence form equations -- and require strong regularity assumptions in both variables. 
	As has been discussed above, \cite{Chen2024} establishes Hölder regularity and a weak Harnack inequality for very general kernels. 
	
	The most relevant result related to \autoref{thm:regularity} is contained in \cite{BassChenRegularity}. The authors consider the operator
	\begin{align}\label{appl:def:bass-chen-generator}
		\mathcal L u(x) 
		\defeq \sum_{k=1}^d c \int_{\R\setminus\{0\}} 
		\big(u(x + a_j(x)t)  - u(x) - 
		t \indic_{B_1}\big(ta_j(x)\big) \nabla u(x) \cdot a_j(x)\big)
		\, |t|^{-1-2s} \, \intd t.
	\end{align}
	The family of measures associated with this operator corresponds to the nonsymmetric family of singular measures from \autoref{def:frontview-concrete}, if $\iota(\intd t) = |t|^{-1-2s} \, \intd t$, $T = \infty$, leaving out the additional weight $|\det D \viewMapKT(x)|$ that is needed for symmetrization.
	It is thus the non-divergence form operator associated with the nonsymmetric family of singular measures $\frontview(x, \intd y)$. 
	In \cite{BassChenRegularity}, Hölder continuity of bounded functions which are harmonic with respect to $\mathcal L$ is established. Moreover, it is shown that a Harnack inequality may fail already in the case where $A \equiv \id$.
	Thus, this result can be seen as the non-divergence form equivalent of \autoref{thm:regularity}. 
	Note that there exists a vast literature on operators in non-divergence form, which is not discussed here. 
	
	\medskip 
	
	Finally, non-translation invariant families of singular measures also appear in the context of regular Dirichlet forms and Markov processes. \autoref{thm:comparability-lip} implies that $\quadForm^\Gamma_\domain$ is a regular Dirichlet form on $H^s(\domain)$ for $s \in (0, \frac 12]$ resp. $H^s_0(\domain)$ for $s \in (\frac 12, 1)$, see \cite[Section 3.1]{HeppDiss}. It thus generates a strong Markov process. 
	Let us quickly recap some literature in this area. 
	
	The stochastic process associated with the family of measures \eqref{def:starrer-fall}
	is sometimes called the rectilinear $\alpha$-stable\footnote{In the context of stochastic processes and stochastic differential equations, the order of differentiability is generally written as $\alpha = 2s \in (0, 2)$.}
	process and has gathered some attention in the literature, see \cite{Xu2013},  \cite{KaKiKu2022}, and  \cite{CHZRectStab25}. 
	%Using its associated Dirichlet form, \cite{Xu2013} and \cite{KaKiKu2022} prove heat kernel bounds. Their results also hold for processes whose Dirichlet forms are merely comparable to that of the rectilinear $\alpha$-stable process. \cite{CHZRectStab25} investigates properties of the transition density function. 
	The rectilinear $\alpha$-stable process can be used as a building block to construct more general jump processes via stochastic differential equations resp. the martingale problem. One may consider a stochastic process $Z_t = (Z_t\upperInd{1}, \dots, Z_t\upperInd{d})$, where $Z_t\upperInd{k}$ are independent one-dimensional symmetric $\alpha$-stable processes. 
	It  drives a stochastic differential equation
	\begin{align}\label{appl:eq:spde-cylindrical-process}
		d X_t = A(Z_{t-}) d Z_t, X_0 = x \in \Rd,
	\end{align}
	where $A: \Rd \to \Rd \times \Rd$ is a matrix-valued function that is bounded and non-degenerate (i.e. $\det A(x) \neq 0$ for all $x \in \Rd$) and satisfies some smoothness assumption. 
	It is well-known that \eqref{appl:eq:spde-cylindrical-process} has a unique and strong solution if $A$ is Lipschitz (see \cite[Theorem 34.7, Corollary 35.3]{MetSemimartingales}). 
	In \cite{BassChenProcess}, it is shown that it has a unique weak solution, even if $A$ is merely continuous. In this case, the infinitesimal generator is given by \eqref{appl:def:bass-chen-generator}.
	Whereas the existence of a process in our case follows from the existence of the regular Dirichlet form in a simple way, the authors of \cite{BassChenProcess} need to invest some effort to solve the martingale problem to deduce the existence of a stochastic process. 
	
	This process and its corresponding stochastic differential equation have been investigated in several other works. Regularity estimates under various smoothness assumptions on $A$ are established in \cite{KulRyz2018}, \cite{ChenEtAl2020}, and \cite{KulEtAl2021}. 
	\cite{PriZa2011} investigates further structural properties of \eqref{appl:eq:spde-cylindrical-process}. 
	
\subsection{Notation and preliminaries}\label{sec:intro:notation}
	We review some notation. 
	We generally work in $\Rd$, where $d$ always denotes the dimension. 
	The unit sphere in $\Rd$ is $\Sd$. 
	The identity matrix in $\Rd$ is $\id \defeq \big(\delta_{ij}\big)_{i, j = 1, \dots, d}$. 
	For $A \in \R^{d\times d}$, the \textit{Frobenius norm} is defined as
	$|| A ||_F \defeq \sqrt{\sum_{i, j = 1}^d |a_{ij}|^2},$
	while the \textit{operator norm} is given by 
	$||A||_{op} \defeq \sup_{x \in \Sd} |Ax|.$
	A ball of radius $R$ around a point $z$ is denoted by $B_R(z)$. In analogy, a square of radius $R$ around $z$ is $Q_R(z) = \bigtimes_{i=1}^d (z_i - R, z_i + R)$. If $z=0$, we often write $B_R$ resp. $Q_R$. Note that, unless explicitly stated otherwise, our balls and cubes are always open. 
	Moreover, recall that $Q_{\frac{R}{\sqrt d}}(z) \subset B_R(z) \subset Q_R(z)$.  If $B = B_R(q)$ is any cube in $\Rd$ and $\lambda > 0$ a factor, then $\lambda B \defeq B_{\lambda R}(q)$ is the ball with radius stretched by $\lambda$. Analogously, we write $\lambda Q \defeq Q_{\lambda R}(q)$ for a cube $Q = Q_R(q)$.
		
	If $F: \Rd \to \R^n$ is a vector-valued function, we write its coordinate functions as $F\upperInd{i}$ for $i \in \{1, \dots, n\}$. If $F$ is differentiable, its Jacobian matrix is 
	\begin{align*}
		\big(D F(x)\big)_{ij} 
		= \frac{\partial F\upperInd{i}}{\partial x_j}
		= \partial_j F\upperInd{i}(x).
	\end{align*}
	For two real numbers $a$ and  $b$, $a \wedge b$ denotes their minimum and $a \vee b$ their maximum.
	Finally, to shorten notation, we sometimes write
	$U(x, y) \defeq \big(u(x) - u(y)\big)^2$
	for points $x, y \in \Rd$. 
	
%	The following lemma is adapted from \cite[proof of Theorem 6.2]{DyKa20}.For the slight generalization here, we introduce the following notation. If $I = \{i_1, \dots, i_k\} \subset \{1, \dots, d\}$ is a set of distinct indices with $k \in \{1, \dots, d-1\}$ and $x \in \Rd$, let 
%	$x^I \defeq (x_{i_1}, \dots, x_{i_k})^T \in \R^k$
%%	\begin{align*}
%%		x^I
%%		\begin{pmatrix}
%%			x_{i_1}\\
%%			\vdots\\
%%			x_{i_k}
%%		\end{pmatrix} 
%%		\in \R^k
%%	\end{align*}
%	and $x^O$ be the $(d-k)$-dimensional vector containing entries with indices in $\{1, \dots, d\} \setminus I$. 
%	If $I = \{k\}$, we write $x^O = x^{\{k\}}$. We refer to \cite[Lemma A.3]{HeppDiss} for a full proof. 
%	\begin{lemma}\label{lemma:barteks-trick}
%		There is a constant $c = c(d) > 0$, such that, for any $s \in (0, 1)$, $k \in \{1, \dots, d-1\}$ with $I$ defined as above
%		\begin{align*}
%			\int_{\R^{d-k}} |x|^{-d-2s} \, \intd x^O
%			\leq c |x^I|^{-k-2s}.
%		\end{align*}
%	\end{lemma}
	
	\medskip
	This paper is structured as follows. 
	\autoref{sec:measures} introduces the families of singular measures $\frontview(x, \intd y)$ and $\fullview (x, \intd y)$ in detail and gives a proof of \autoref{thm:properties-for-regularity-theory}. 
	The last three sections deal with the proof of \autoref{thm:comparability-lip}. We prove the lower bound in \autoref{sec:connectedness} and \autoref{sec:coerc}. \autoref{sec:connectedness} contains a connectedness result for Lipschitz configurations, which is used to prove coercivity in \autoref{sec:coerc}. The upper bound is treated in the final \autoref{sec:upper}. 	
	
\subsection{Acknowledgments}
	We thank Moritz Kassmann for introducing us to this problem and for many insightful discussions.
	Kai-Uwe Bux thanks the Deutsche Forschungsgemeinschaft (DFG, German Research Foundation) for support via the CRC TRR 358 – Project-ID 491392403. 
	Solveig Hepp gratefully acknowledges financial support by the German Research Foundation (GRK 2235 - 282638148). 

\subsection{AI Statement}
	No AI has been used for this article, neither in the generation of the mathematical statements or proofs, nor in the writing. 
	The results presented are part of the PhD project of Solveig Hepp, see \cite{HeppDiss}.

\section{Space-dependent singular measures}\label{sec:measures}
	This section introduces the main objects of this paper: configurations and the families of singular measures they generate.  
	In \autoref{sec:measures:configs}, we discuss the properties \eqref{ass:lip}, \eqref{ass:det}, and \eqref{ass:w} in more detail. 
	\autoref{sec:measures:nonsymmetric-family} establishes some properties of families of measures like $\frontview(x, \intd y)$ given by \autoref{def:frontview-concrete}. In particular, we show that they are intrinsically nonsymmetric.
	For this reason, new, symmetric families of measures $\fullview(x, \intd y)$ are constructed in \autoref{sec:measures:symmetric-family}, where we also prove \autoref{thm:properties-for-regularity-theory}.
	
\subsection{Configurations}\label{sec:measures:configs}
	There are broadly three types of interesting properties a configuration may possess: the properties of each subconfiguration on a small scale (e.g. \eqref{ass:lip}), the properties of each subconfiguration on a large scale (e.g. \eqref{ass:w}), and the way the subconfigurations relate to each other (e.g. \eqref{ass:det}). 
	
	The \textit{small-scale behavior} of a configuration is determined by the regularity of its subconfigurations. 
	Therefore, unless otherwise noted, from now on $\Gamma: \Rd \to \Sd$ will always signify a subconfiguration. 
	The importance of a subconfiguration's regularity lies in the fact that it governs the regularity of the view map. 
	In this context, the goal is to control the pre-image of each point $x$ under $\viewMap_t$, at least in some neighborhood of $x$. 
	While invertibility of the view map on the whole space will, in general, only be true for constant configurations, invertibility on small scales will hold if $\Gamma$ is Lipschitz.
	
	Let us prove a basic result about view maps of Lipschitz configurations. 
	\begin{lemma}\label{measures:lemma:injectivity-basic-for-lipschitz-configs}
		For each $\constL > 0$, there is a radius $R = R(\constL) > 0$, such that, for any Lipschitz continuous subconfiguration $\Gamma: \Rd \to \Sd$ with Lipschitz constant $\constL$,
		and any $t \in \R$ and $x, x' \in \Rd$ with $|t|, |x - x'| < R$, the following is true: If $\viewMap_t(x) = \viewMap_{t}(x')$, then $x = x'$. 
	\end{lemma}
	\begin{proof}
		If $\Gamma(x) = \Gamma(x')$, the assertion holds trivially, independent of any bounds on $|t|$ or $|x-x'|$. So assume $\Gamma(x) \neq \Gamma(x')$. 
		Let $R_1 = R_1(\constL) > 0$ be such that 
		\begin{align*}%\label{meas:eq:injectivity-basic-choice-radius-1}
			\arccos\Big( 1 - \frac{\constL^2 R_1^2}{2}\Big) < \frac \pi 4
			\text{ and }
			R_1 < \frac{\sin\big(\frac \pi 4\big)} {\constL}.
		\end{align*}
		Moreover, as
		\begin{align*}
			\frac{a_n}{\sin\Big(\arccos \Big(1 - \frac{1}{2} c^2 a_n^2 \Big)\Big)}
			\searrow \frac{1}{c}
			\text{ for any sequence }
			(a_n)_{n \in \N} \searrow 0,
		\end{align*}
		we find $R \leq R_1$, such that:
		\begin{align*}%\label{meas:eq:injectivity-basic-choice-radius-2}
			\forall \; a \in (0, R) \; : \; 
			\frac{a}{\arccos \big( 1 - \frac{1}{2} \constL^2 a^2\big)} > \frac{1}{\constL}
		\end{align*} 
		Assuming for the sake of contradiction that $|t|, |x-x'| < R$, $\viewMap_t(x) = \viewMap_{t}(x')$, but $x \neq x'$, the assertion follows from the choice of $R$ and the law of sines applied to the isosceles triangle with edges $x, x'$, and $\viewMap_t(x) = \viewMap_t(x')$. 
	\end{proof}
	
	Next, we derive some useful properties of the view map of a Lipschitz subconfiguration. 
	\begin{lemma}\label{measures:lemma:properties-of-det-view-map}
		Let $\Gamma: \Rd \to \Sd$ be a subconfiguration that satisfies \eqref{ass:lip} with Lipschitz constant $\constL > 0$. Then $\viewMap_t: \Rd \to \Rd$ is Lipschitz continuous for each $t \in \R$.
		There also holds: 
		\begin{enumerate}[(a)]
			\item There is a constant $c = c(d, \constL) > 0$, such that, for every $x \in \Rd$ where $D \viewMap_t(x)$ exists,
			$|\det D \viewMap_t(x)| \leq c \,\big(1 \vee |t|^d\big).$
			\item For each $\delta > 0$, there is a radius $R_0 = R_0(d, \constL, \delta) > 0$, such that, for every $|t| < R_0$ and every $x \in \Rd$ where $D \viewMap_t(x)$ exists,
			$||D \viewMap_t(x) - \id ||_F < \delta.$ 
			\item For each $\kappa \in (0, 1)$, there is a radius $R_1 = R_1(d, \constL, \kappa) > 0$, such that, for every $|t| < R_1$ and every $x \in \Rd$ where $D \viewMap_t(x)$ exists,
			$|\det D \viewMap_t(x)|  > \kappa.$
			\item There is a radius $R_2 = R_2 (\constL) > 0$, such that, for any $z \in \Rd$, and $|t| < R_2$, $\viewMap_t: B_{R_2}(z) \to \Rd$ is injective. 
		\end{enumerate}
	\end{lemma}
	\begin{proof} 
		The Lipschitz continuity of $\viewMap_t$ follows from \eqref{ass:lip}. 
		Due to Rademacher's theorem (see e.g. \cite[Theorem 5.8.6]{Evans10}), the partial derivatives of a Lipschitz continuous subconfiguration $\Gamma$ exist in almost every point and are bounded by $\constL$. Let $x \in \Rd$ be such a point. 
		Observe that $\big(D \viewMap_t(x)\big)_{ij} = \big(\delta_{ij} + t \partial_j \Gamma\upperInd{i}(x)\big)$. The boundedness of the partial derivatives of $\Gamma$ gives (a).
		Since 
		\begin{align*}
			||D \viewMap_t(x) - \id ||_F
			= \sqrt{\sum_{i, j = 1}^d \big(\delta_{ij} + t \partial_j \Gamma\upperInd{i}(x) - \delta_{ij}\big)^2}
			\leq d R_0 \constL, 
		\end{align*}
		we may choose $R_0$ depending only on $d, \constL,$ and $\delta$ to obtain (b), from which (c) follows together with the continuity of the determinant. 
		Finally, let $R$ be the radius from \autoref{measures:lemma:injectivity-basic-for-lipschitz-configs} and set $R_2 \defeq \frac{R}{2}$. Now assume $B\subset \Rd$ is any ball of radius at most $R_2$ and $|t| < R_2$. If $x, x' \in B$ with $\viewMap_t(x) = \viewMap_t(x')$, then \autoref{measures:lemma:injectivity-basic-for-lipschitz-configs} implies $x =x'$, since $|x-x'| < 2 R_2 < R$. 
		Thus $\viewMap_t$ is injective on $B$, i.e. (d) holds. 
	\end{proof}
	
	\medskip
	
	Let us now discuss the \textit{large-scale properties} that $\Gamma$ bestows upon $\viewMap_t$ for large $|t|$.
	Even for the most well-behaved subconfigurations, $\viewMap_t(x)$ and $\viewMap_t(y)$ may be arbitrarily far apart from each other for $|t| \to \infty$, no matter how close $x$ and $y$ are. 
	This might cause a situation, where a non-negligible set of points is seen at a non-negligible set of distances by infinitely many other non-negligible sets of points. In such a scenario, certain sets of points can obtain too large a weight, which might cause problems when proving the upper bound of \eqref{eq:comparability} later on. 
	
	Property \eqref{ass:w} prevents this kind of phenomenon. 
	It is independent of the regularity property \eqref{ass:lip}, as can be seen from the fact that configurations with finite image always satisfy \eqref{ass:w} with $M = \#\{\Gamma(\Rd)\}$. 
	On the other hand, in the following example, a subconfiguration is constructed that can be arbitrarily regular but does not satisfy \eqref{ass:w}. 
	\begin{example}
		\label{measures:ex:wurst-example}
		Let $r \in (0, 1)$ and consider $B_r \subset \R^2$. 
		Let $\rho_n \defeq 2^{-n}r$ for each $n \in \N$ be a sequence of auxiliary radii. 
		Furthermore, let $\thet_n \defeq \arctan \big(\frac{\rho_n}{2r}\big) = \arctan (2^{-n-1})$. 
		If $\phi \in (0, \pi)$ is an angle, we write 
		$\vect v(\phi) \defeq \begin{pmatrix}
				\cos(\phi)\\
				\sin(\phi) 
			\end{pmatrix}
			\in \mathbb S^1,$
		for the unit vector defined by this angle, and set 
		$\phi_0 \defeq \thet_0
			\text{ and }
			\phi_n \defeq \sum_{i=1}^{n-1} 3 \thet_i + \thet_n.$
		For a radius $\rho$, an angle $\phi$, and a point $x \in \R^2$, define the tube
		\begin{align*}
			\wurst (\rho, \phi, x)
			\defeq 
			\{ B_\rho(x) + t \vect v(\phi) \mid t \in (2 r, 3 r)\}
		\end{align*}
		and set $\wurst_n \defeq \wurst(\rho_n, \phi_n, 0)$. 
		Observe that $\dist (\wurst_n, \wurst_m) > 0$ if $n \neq m$. 
				
		Since $\sum_{n=1}^\infty 3 \thet_n = \sum_{n=1}^\infty 3 \arctan(2^{-n-1}) < \pi$, the tubes $\wurst_n$ all live on the upper half-plane and in $B_{4r}$. 
		Define the subconfiguration on each tube as
		\begin{align*}
			\Gamma(x) \defeq - \vect v(\phi_n) 
			\text{ if } x \in \wurst_n.
		\end{align*}
		Now interpolate between the tubes in any way you want to produce a configuration of any smoothness or even a configuration with countable image. 
		
		For any $y \in B_{\rho_n}$ and $t \in (2r, 3r)$, we have 
		\begin{align*}
			\haus^0 \big(\{x \in \Rd \mid \viewMap_t(x) = y\}\big)
			\geq n,
		\end{align*}
		since $B_{\rho_n} + t \vect v(\phi_n) \subset \wurst_k$ for any $k \in \{1, \dots, n\}$. 
		Thus, for $t \in (2r, 3r)$, 
		\begin{align*}
			\essup_{y \in B_{4r}}\;\haus^0 \big(\{x \in \Rd \mid \viewMap_t(x) = y\}\big) = \infty.
		\end{align*}
		However, \eqref{ass:w} would require $\haus^0 \big(\{x \in \Rd \mid \viewMap_t(x) = y\}\big)$ to be bounded for $y \in B_{4r}$ and $t \in (2r, 3r)$. 
		Hence, $\Gamma$ does not satisfy \eqref{ass:w}. 
	\end{example}
	
	\medskip
	
	For the properties concerning the \textit{relationship of the subconfigurations to each other}, let us assume that $\Gamma: \Rd \to \big({\Sd}\big)^d$ is a full configuration once more. 
	We will call $\Gamma_k(x)$ the $k$-th \textit{type vector} of $x$. 
	The \textit{type vector matrix} of $x$ is the matrix that contains the type vectors of $x$ as its columns, i.e. 
	\begin{align*}%\label{def:rot-m}
		\RotM_x \defeq 
		\begin{pmatrix}
			|		& 		& | \\
			\Gamma_1(x)	& \cdots& \Gamma_d(x)\\
			|	&	&|
		\end{pmatrix}.
	\end{align*}
	Observe that $\RotM_x$ will always be invertible. Its inverse $\RotM_x\inv$ is the matrix changing the standard basis $(\vect e_1, \dots, \vect e_d)$ to $\Gamma(x)$, i.e. $\RotM_x(\vect e_k) = \Gamma_k(x)$ and, for any vector $\vect v \in \Rd$, $\vect v = \sum_{k=1}^d \RotM_x\inv(\vect v)_k \Gamma_k(x)$.	
	
 	The fact that $\RotM_x$ is a matrix whose columns consist of unit vectors implies $|\det \RotM_x| \leq 1$ and $||\RotM_x||_{op} \leq ||\RotM_x||_F \leq 2^\frac{d}{2} d$, hence 
	$||\adj \RotM_x||_F \leq d^\frac{d-2}{2} (2^\frac{d}{2}d)^{d-1} = c(d)$ (see \cite{Richter54}). The precise value of the constant is not important, since it only depends on the dimension $d$. 
	The nondegeneracy assumption \eqref{ass:det} stated in the introduction essentially serves to derive a bound on the Frobenius norm of the inverse of $\RotM_x$ as well. 
	It implies, together with the estimates above, 
	\begin{align}\label{meas:cc:config-det:eq:norm-of-inverse-bound}
		||\RotM_x \inv||_{op} 
		\leq ||\RotM_x \inv||_F
		= |\det \RotM_x|\inv ||\adj \RotM_x||_F
		\leq \constDet\inv c(d)\tag{Det2}
	\end{align}
	for any type vector matrix $\RotM_x$. 
	Thus, the uniform lower bound on the determinant of $\RotM_x$ implies a uniform upper bound on $||\RotM_x\inv||_F$. 
	
\subsection{The nonsymmetric family of singular measures}\label{sec:measures:nonsymmetric-family}
	Let us slightly generalize \autoref{def:frontview-concrete} of $\frontview(x, \intd y)$ given in the introduction. 
	We define a \textit{family of set functions} with cutoff point $T \in (0, \infty]$ generated by a configuration $\Gamma: \Rd \to \big({\Sd}\big)^d$, a one-dimensional $\sigma$-additive Lebesgue measure $\iota(\intd t)$, and a family of functions $F_k: \R^{d+1} \to [0, \infty]$ by
	\begin{align}\label{measures:def:frontview-generalized-definition}
		\forall \; A \in \borel(\Rd), \; x \in \Rd \; : \; 
		\nu(x, A) 
		\defeq \sum_{k=1}^d \int_{-T}^T 
		\indic_A \big(x + t \Gamma_k(x)\big) F_k(x, t) \, \iota(\intd t). 
	\end{align}
	When establishing properties of these set functions, it is once again sufficient to only consider subconfigurations.
	The following technical observation follows from basic properties of measurable functions. 
	\begin{lemma}[Measurability of the view map]\label{measures:lemma:view-map-is-measurable}
		If $\Gamma: \Rd \to \Sd$ is $\borel(\Rd)-\borel(\Sd)$-measurable, then the map $V: 	(x, t) \mapsto x + t \Gamma(x)$
		is $\borel(\R^{d+1})-\borel(\Rd)$-measurable. 
	\end{lemma}
	We may now derive properties of measures $\nu(x, \intd y)$ from the properties of the objects that generate them.
	\begin{lemma}\label{measures:lemma:nonsymmetric-measure-well-defined}
		Let $\Gamma: \Rd \to \Sd$ be a subconfiguration, 
		$\iota(\intd t)$ a one-dimensional $\sigma$-additive Lebesgue measure, 
		$F(x, t)$ a $\leb(\R^{d+1}) - \borel([0, \infty])$-measurable, non-negative function, and $\nu(x, \intd y)$ the family of set functions with cutoff point $T \in (0, \infty]$ generated by them according to \eqref{measures:def:frontview-generalized-definition}. 
		\begin{enumerate}[(a)]
			\item For any $x \in \Rd$, $\nu(x, \fctempty)$ is a Borel measure on $\Rd$. For any $A \in \borel(\Rd)$, $\nu(\fctempty, A)$ is $\leb(\Rd)-\borel([0, \infty])$-measurable. 
			\item If, additionally, $\iota(\intd t)$ satisfies \eqref{def:iota-comparable} for $s \in (0, 1)$ 
			and $F$ is uniformly bounded, $\nu(x, \intd y)$ satisfies \eqref{def:decay} for $s$, as well as \eqref{def:levy}. 
		\end{enumerate}
	\end{lemma}
	
	\begin{proof}
		The second assertion of (a) can be verified directly using basic properties of the characteristic function $\indic_A$. 
		The first assertion follows from \autoref{measures:lemma:view-map-is-measurable} together with Tonelli's theorem. 
		The proof of (b) is a simple calculation. Assume $\iota(\intd t)$ satisfies \eqref{def:iota-comparable} with $\constI$ and suppose $\sup |F(x, t)|$ $< c_0 < \infty$. If $x_0 + t \Gamma(x_0) \in B_r(x_0)^c$, $|t| > r$. Thus
		\begin{align*}
			\nu\big(x_0, B_r(x_0)^c\big) 
			&= \int_{\{r < |t| < T\}} \indic_{B_r(x_0)^c}\big(x_0 + t \Gamma(x_0)\big)
			F(x_0, t) \, \iota(\intd t)\\
			&\leq \; c(\constI, c_0) \; \int_r^\infty t^{-1-2s} \,\intd t 
			= \frac{c(\constI, c_0)}{2s} r^{-2s}.
		\end{align*}
		The same kind of calculation can be used to deduce \eqref{def:levy}. 
	\end{proof}
	
	\begin{remark}\label{meas:rem:measurability-of-the-det-function}
		If $\Gamma: \Rd \to \Sd$ is a Lipschitz continuous subconfiguration, the map $(x, t) \mapsto {|\det D \viewMap_t(x)|}$
		is $\leb(\R^{d+1})$-measurable. 
		Thus, \autoref{measures:lemma:nonsymmetric-measure-well-defined} (a) implies the well-definedness of the family of measures $\frontview(x, \intd y)$ given in \autoref{def:frontview-concrete}.
		If the cutoff point is chosen in $(0, \infty)$, $|\det D \viewMap_t(x)|$ will be uniformly bounded due to \autoref{measures:lemma:properties-of-det-view-map} (a), which, together with \autoref{measures:lemma:nonsymmetric-measure-well-defined} (b), gives \eqref{def:decay} and \eqref{def:levy} for $\frontview(x, \intd y)$ as well. 
	\end{remark}
	
	No matter the choices of $\Gamma$, $\iota(\intd t)$, and $F$, however, the family of measures $\nu(x, \intd y)$ will, in general, not be symmetric, as the following example demonstrates. 
	\begin{example}[Nonsymmetric family of measures]
		\label{meas:ex:nonsymmetric-measure}
		Let $A = (-4, -1) \times (-3, -2)$, $B = (1, 4) \times (2, 3)$. Define a configuration $\Gamma: \R^2 \to (\mathbb S^1)^2$ by setting $\Gamma_1 \equiv \vect e_1$ and 
		\begin{align*}
			\Gamma_2 (x)
			\defeq 
			\begin{cases}
				\vect e_2 	&\text{on }A\\
				\begin{pmatrix}
					\frac{1}{\sqrt 2}\\
					\frac{1}{\sqrt 2}
				\end{pmatrix}
				&\text{on }B. 
			\end{cases}
		\end{align*}
		%See also Figure 4.2. 
		It is possible to extend $\Gamma_2$ to all of $\Rd$ in such a way that the resulting configuration has any regularity, or, alternately, a finite image, and \eqref{ass:w} and \eqref{ass:det} hold. 
		However, for all but a few very special choices of $F$ and $\iota(\intd t)$, 
		\begin{align*}
			\int_A \int_B \nu (x, \intd y)\, \intd x = 0 \text{ but } \int_B \int_A \nu(x, \intd y) \,\intd x > 0.
		\end{align*}
	\end{example}
%	\begin{figure}[h!]
%		\centering
%		\includegraphics[width=\textwidth/2]{../diss/pics/nonsymmetric.png}
%		\caption[Short caption]{A configuration for which $\nu(x, \intd y)$ is not symmetric.}
%	\end{figure}

\subsection{A symmetric family of singular measures}\label{sec:measures:symmetric-family}
	We will now construct the symmetric family of measures given in \autoref{thm:properties-for-regularity-theory}, which implies the regularity results of \autoref{thm:regularity}. 
	We have seen that $\nu(x, \intd y)$ is not symmetric in the sense of \eqref{def:symmetry}. 
	But it gives rise to a measure $n$ on the product space that satisfies
	\begin{align}\label{def:mGam}
		\forall \; A, B \in \borel(\Rd) \; : \; 
		A \cap B = \emptyset \; 
		\Rightarrow \; 
		n(A, B) = \int_A \nu(x, B) \, \intd x\tag{$P_1$}
	\end{align}
	and can be symmetrized to yield a measure $m$ with
	\begin{align}
		\label{def:mGamSym}
		\forall \; A, B \in \borel(\Rd)\; : \;
		A \cap B = \emptyset \; \Rightarrow \;
		m(A, B) &\defeq \tfrac{1}{2}\Big(\int_A \nu(x, B) \, \intd x + \int_B \nu(x, A) \, \intd x\Big).\tag{$P_2$}
	\end{align}
	Now $n$ corresponds to the family of measures $\nu(x, \intd y)$ in the sense that, for each Borel set $B$, $\nu(\fctempty, B)$ acts as the density of $n(\fctempty, B)$ with respect to the Lebesgue measure. 
	Analogously, if there is a family of measures $\mu(x, \intd y)$ that satisfies
	\begin{align}\label{eq:density}
		\forall \; A, B \in \borel(\Rd) \; : \; 
		A \cap  B = \emptyset
		\; \Rightarrow \; 
		m(A, B) 
		= \int_A \mu(x, B) \, \intd x,\tag{Dens}
	\end{align}
	we will say that the family of measures $\nu(x, \intd y)$ is \textit{symmetrizable} and call $\mu(x, \intd y)$ the \textit{family of symmetric singular measures} generated by $\Gamma$ and $\iota$ and cut off at $T$. Clearly, if $\mu(x, \intd y)$ satisfies \eqref{eq:density}, it also satisfies \eqref{def:symmetry}. 
	
	\medskip
	
	Let us begin by providing an example of a configuration generating a family of measures $\nu(x, \intd y)$ that is \textit{not} symmetrizable. 
	The main idea is this. If \eqref{eq:density} is to hold, then the measure $m(\fctempty, B)$ must have a density $f_B(x) = \mu(x, B)$ with respect to the Lebesgue measure for every $B \in \borel(\Rd)$. 
	Thus, each $m(\fctempty, B)$ must be absolutely continuous with respect to the Lebesgue measure. 
	Conversely, if some $m(\fctempty, B)$ is \textit{not} absolutely continuous for some $B$, there cannot be a density $f_B$ and, consequently, there cannot be a well-defined family of measures $\mu(x, \intd y)$ satisfying \eqref{eq:density}. 
	
	A Nikod\'ym set, introduced by Otton Nikod\'ym in his \cite{Nikodym1927}, is a Borel set of full measure in some bounded set $Q$, with a corresponding null set that contains lines such that each point in the Nikod\'ym set lies on at least one such line. 
	These lines can even be chosen in a measurable way, as the result \cite[Theorem 8.5.5]{Guz1981} demonstrates: Let $Q = (-1, 1)^2$ be the unit square in $\R^2$. There exist Borel sets $B \subset Q$ with $|B| = |Q|$, and $K \subset \R^2$ with $|K| = 0$, such that, for any $x \in B$, there is a line $L_x = x + \langle \vect v(x) \rangle$ passing through $x$, such that $L_x  \subset K \cup \{x\}$. Moreover, the \textit{accessibility function} $\vect v: B \to \mathbb S^1$, which assigns to each point the direction of its line, is Borel measurable.
	The set $B$ is thus a Nikod\'ym set with measurable accessibility function. 
	Using this result, we construct a family of measures that cannot be symmetrized. 
	
	\begin{example}[Non-symmetrizable family of measures]
		\label{meas:ex:non-symmetrizable-family-of-measures}
			Let $Q$, $B$, and $K$ be as above. We may take $B$ and $K$ to be disjoint. Set $N \defeq Q \setminus B \supset Q \cap K$. Since $|B| = |Q|$, $|N| = 0$. Define 
		\begin{align}\label{meas:eq:nikodym-set-configuration}
			\Gamma (x) = 
			\begin{cases}
				\vect v(x) 	& \text{if } x \in B\\
				\vect e_1 	& \text{if } x \in N.
			\end{cases}
		\end{align}
		Then $\Gamma$ may be extended to a measurable function on all of $\R^2$ and is thereby a valid subconfiguration. For simplicity, let $\nu(x, \intd y)$ be the family of measures generated by $\Gamma$, the one-dimensional Lebesgue measure, and $F \equiv 1$. It is clear that the same construction works for the family $\frontview(x, \intd y)$ of \autoref{def:frontview-concrete}. 
		
		For any $x \in B$, $x + t \Gamma(x) = x + t \vect v (x) \in N$ if $t \neq 0$. In particular, if $x \in (-\frac{1}{2}, \frac{1}{2})^2$, $B_{\frac{1}{2}}(x) \subset Q$, and thus 
		\begin{align*}
			\int_\R \indic_N \big(x + t \Gamma(x)\big) \, \intd t 
			\geq 1 > 0. 
		\end{align*}
		As $|B| = |Q|$, this implies
		\begin{align*}
			\int_B \nu(x, N) \, \intd x  
			= \int_B \int_\R \indic_N \big(x + t \Gamma(x)\big) \, \intd t \intd x 
			\geq |B \cap (-\tfrac{1}{2}, \tfrac{1}{2})^2| > 0.
		\end{align*}
		Thus, $m(N, B) > 0$ but $|N| = 0$, so $m(\fctempty, B)$ is not absolutely continuous with respect to the Lebesgue measure and therefore cannot have a density either. 
	\end{example}
	
	Now we construct a symmetric family of measures. 
	The different values for $\int_A\nu(x, B) \, \intd x$ and $\int_B\nu(x, A)\, \intd x$ in \autoref{meas:ex:nonsymmetric-measure} are due to the fact that $B$ sees $A$, but $A$ cannot see $B$. 
	Hence, in order to symmetrize $\nu(x, \intd y)$, one would like to add a term to $\int_A \nu(x, B) \, \intd x$ that takes into account everything $B$ sees as well. In fact, this is exactly what the second summand of the symmetric product measure $m$ characterized by \eqref{def:mGamSym} does. However, in order to get a family of measures of the form $\mu(x, \intd y)$, we need the second term to be an integral over $A$, not over $B$. In other words, we want to write something like  
	\begin{align*}
		m(A,  B) =
		\int_A \nu(x, B) \, \intd x 
		+ \int_A \int_\R G_B(x, t) \,\iota(\intd t) \intd x,
	\end{align*}
	where $G_B(x, t)$ is some measurable function depending on the set $B$. 
	
	A natural candidate for $G_B$ is the multiplicity function of each view map $\viewMap_t$,
	\begin{align*}
		\Rd \times \R &\to \N \cup \{\infty\}\\
		(x, t) &\mapsto \haus^0\big( \{y \in B \mid \viewMap_t(y) = x\}\big)
	\end{align*}
	which counts the number of points $y \in B$ that see $x$ at distance $t$. When integrating this function over $A$ and $\R$, the distances are measured using $\iota(\intd t)$ and weighted with the number of preimages of each $x \in A$. This is exactly the behavior we would expect of the second summand in the symmetrized measure $m$. 
	We therefore define: 
	\begin{definition}\label{meas:def:rearview-measure}
		Let $\Gamma: \Rd \to (\Sd)^d$ be a configuration, 
		$\iota: \leb(\R) \to [0, \infty]$ a one-dimensional $\sigma$-additive Lebesgue measure, 
		and $T \in (0, \infty]$. 
		The \textit{family of rear view measures} with cutoff point $T$ generated by $\Gamma$ and $\iota$ is given by 
		\begin{align*}
			\rho(x, \fctempty): \borel(\Rd) &\to [0, \infty]\\
			A & \mapsto  
			\sum_{k=1}^d \int_{-T}^T
			\haus^0\big(\{ y \in A \mid \viewMapKT(y) = x\}\big) \, \iota(\intd t).
		\end{align*}
	\end{definition}		
	As before, we have to verify the well-definedness of $\rho(x, \intd y)$. 
	
	In order to prove \eqref{def:well-definedness-measures}, we need the multiplicity function to be $\leb(\R^{d+1})-\borel([0, \infty])$-measurable for every Borel set $A$. This will probably not be the case in general. We do, however, have the following result. 
	\begin{lemma}\label{meas:lemma:criterion-for-mb}
		Let $\Gamma: \Rd \to \Sd$ be a continuous configuration, $I \subset \R$ an interval, and $A \in \borel(\Rd)$. Assume there are pairwise disjoint sets $(M_n)_{n \in \N} \subset \borel(\Rd)$ with $\bigcup_{n \in \N} M_n = A$, such that 
		\begin{align*}
			\forall \; n \in \N, \text{ a.e. } t \in I \; : \; 
			\viewMap_t: M_n \to \Rd \text{ is injective}. 
		\end{align*}
		Then the multiplicity map $(x, t) \mapsto \haus^0 \big(\{y \in A \mid \viewMapKT(y) = x\}\big)$ is $\leb(\Rd \times I)-\haus^0(\N \cup \{\infty\})$-measurable. 
	\end{lemma}
	\begin{proof}
		Fix one $M \in \borel(\Rd)$, and assume $\viewMap_t : M \to \Rd$ is injective for any $t \in I \setminus B$, with $|B| = 0$. 
		Observe that the multiplicity map for $M$ equals $\indic_{\bigcup_{t \in I}\viewMap_t(M) \times \{t\}} (y, t)$ for $(y, t) \in \Rd \times (I \setminus B)$, i.e. up to a $\leb(\R^{d+1})$-null set.
		Since that map $H: (x, t, s) \mapsto \big(\viewMap_t(x), s\big)$ is continuous, observing $\bigcup_{t \in I}\viewMap_t(M) \times \{t\} = H \big(M\times \diag(I \times I)\big)$ finishes the proof. 
	\end{proof}
	
	Now we can show that $\rearview(x, \intd y)$ is well-defined. 
	\begin{lemma}
		\label{measures:lemma:rearview-measure-is-well-defined}
		Let $\Gamma: \Rd \to \Sd$ be a subconfiguration and 
		$\iota(\intd t)$ a one-dimensional $\sigma$-additive Lebesgue measure. 
		\begin{enumerate}[(a)]
			\item For any $x \in \Rd$ and cutoff point $T \in (0, \infty]$, the set function 
			$\rho (x, \fctempty): \borel(\Rd) \to [0, \infty]$ generated by $\Gamma$ and  $\iota$ is a Borel measure.
			\item If $\Gamma$ additionally satisfies \eqref{ass:lip} with Lipschitz constant $\constL > 0$, there is a cutoff point $T_0 = T_0(d, \constL) > 0$, such that, for any $A \in \borel(\Rd)$, $(x, t) \mapsto \haus^0\big(\{y \in A \mid \viewMap_t(y) = x\}\big)$ is $\leb (\Rd \times (-T, T))$-measurable, and 
			$\rho(\fctempty, A)$ is $\leb(\Rd)-\borel([0, \infty])$-measurable, if $\rearview(x, \intd y)$ is generated by $\Gamma$ and $\iota(\intd t)$ with any cutoff point $T \in (0, T_0)$.
		\end{enumerate}
	\end{lemma}
	\begin{proof}
		By \autoref{measures:lemma:view-map-is-measurable}, we have for any $x \in \Rd$ and $A \in \borel(\Rd)$, 
		\begin{align*}
			M \defeq \viewMap\inv\big( \{x\}\big) \cap (A \times \R) 
			\in \borel(\Rd \times \R).
		\end{align*}
		Since $\{y \in A \mid \viewMap_t(y) = x\}$ is the $t$-slice through $M$, the measurability of the multiplicity map follows. This gives (a). 
		For (b), we only need to construct an upper bound $T_0$, such that $(x, t) \mapsto \haus^0 \big(\{y \in A \mid \viewMap_t(y) = x\}\big)$
		is $\leb\big(\Rd \times (-T, T)\big)$-measurable for any $A \in \borel(\Rd)$ and $T < T_0$. 
		Due to \autoref{measures:lemma:properties-of-det-view-map} (d), there is a radius $R_0 = R_0(d, \constL) > 0$, such that $\viewMap_t: B_{R_0}(z) \to \Rd$ is injective if $|t| < R_0$ for any $z \in \Rd$. Since $Q_{\frac{R_0}{2 \sqrt d}}(z) \subset B_{R_0}(z)$, we may cover $\Rd$ with countably many half-open disjoint cubes of radius $\frac{R_0}{2 \sqrt d}$. This yields the assumption of \autoref{meas:lemma:criterion-for-mb} for any $A \in \borel(\Rd)$, from which the measurability of any view map follows. 
	\end{proof}
	
	In order to prove \eqref{def:decay} and \eqref{def:levy} for the family of rear view measures $\rho(x, \intd y)$, we not only need to assume that $\iota(\intd t)$ behaves like $|t|^{-1-2s}\intd t$, but also that $\Gamma$ satisfies \eqref{ass:w}. Let us first show that \eqref{ass:w} is a necessary assumption for \eqref{def:decay}.
	
	\begin{example}\label{measures:ex:decay-counterexample}
		We build off of \autoref{measures:ex:wurst-example} to construct a configuration that generates a family of rear view measures for which the decay condition \eqref{def:decay} is false for a sequence $(r_N)_{N \in \N}$ of radii with $r_N \to 0$. Specifically, we will construct a subconfiguration $\Gamma: \Rd \to \Sd$ with the following property: If $\rho(x, \intd y)$ is the family of rear view measures generated by $\Gamma$ and $\iota(\intd t) = |t|^{-1-2s}\, \intd t$, we find a constant $c = c(s) > 0$, such that, for each $N, n \in \N$, there is a set $M^{(N)}(n) \subset \domain$ with $|M\upperInd{N}(n)| > 0$ and 
		\begin{align*}
			\forall \; x_0 \in M^N(n) \; : \; \rho\big(x_0, B_{r_N}(x_0)^c\big) 
			> c(s) \, n r_N^{-2s}.
		\end{align*}
		Since $c$ is independent of $N$ and $n$, this implies the negation of \eqref{def:decay}. 
		Note that it suffices to show such an assertion for a subconfiguration since the measure generated by a full configuration dominates the measure generated by one of its subconfigurations.  
		
		We construct a subconfiguration $\Gamma$ that contains countably many shrinking copies of \autoref{measures:ex:wurst-example} and interpolate between the copies. 
		So for $N \in \N$, set $r_N \defeq 2^{-N}$ and consider $B_N \defeq B_{4r_N}(q_N)$ for some $q_N \in \R^2$ such that $\dist(B_N, B_M) > 0$ if $N \neq M$. Since 
		\begin{align*}
			\sum_{n=1}^\infty |B_{8r_N}(q_N)| 
			= |B_1| \sum_{N=1}^\infty (8 r_N)^2 
			= 64 |B_1| \sum_{N=1}^\infty 2^{-2N} 
			< \infty,
		\end{align*}
		we find a bounded domain $\domain$ that contains every $B_N$ in such a way that allows them to have positive distance. 
		
		Let $\rhoN_n \defeq 2^{-N-n}$, $\thetN_n \defeq \arctan \big(\rhoN_n/(2r_N)\big) = \arctan (2^{-n-1})$ and $\phiN_0 \defeq \thetN_0$, $\phiN_n \defeq \sum_{i=1}^{n-1} 3 \thetN_i + \thetN_n.$
		Define $\wurstN_n \defeq \wurstN(\rhoN_n, \phiN_n, q_N)$. 
		These tubes all have positive distance and are contained in the upper half-plane through $q_N$. More precisely $\wurstN_n \subset B_N \setminus B_{r_N}(q_N) \cap \R^2_+$. 
		Define the configuration on each tube as
		\begin{align*}
			\Gamma(x) \defeq - \vect v(\phiN_n) 
			\text{ if } x \in \wurstN_n,
		\end{align*}
		and interpolate between the tubes in whatever way you see fit. 
		
		By \autoref{measures:ex:wurst-example}, we know that $\haus^0\big(\{x \in \Rd \mid \viewMap_t(x) = y\}\big) \geq n$ for any $y \in B_{\rhoN_n}(q_N)$. We obtain for any $r_N > 0, n \in \N$:
		\begin{align*}
			\essup_{x_0 \in \Rd} \rho\big(x_0, B_r(x_0)^c\big)
			&\geq \essup_{x_0 \in B_{\rhoN_n}(q_N)} 
			\int_{2r_N}^{3r_N} \haus^0\big(\{ y \in \Rd \mid \viewMap_t(y) = x_0\}\big) |t|^{-1-2s} \, \intd t\\
			&\geq 
			n \int_{2r_N}^{3r_N} t^{-1-2s} \, \intd t
			\geq c(s) \, n \; r_N^{-2s}.
		\end{align*}
		This implies the negation of \eqref{def:decay}. 
	\end{example}
	
	We now show \eqref{def:decay} and \eqref{def:levy} for $\rho(x, \intd y )$ under the assumption \eqref{ass:w}. 
	\begin{lemma}\label{measures:lemma:rearview-decay-levy}
		Let $\Gamma: \Rd \to \Sd$ be a subconfiguration that satisfies \eqref{ass:w},
		$\iota(\intd t)$ a one-dimensional Lebesgue measure that satisfies \eqref{def:iota-comparable} for $s \in (0, 1)$, and $T\in (0, \infty)$. Any family of rear view measures generated by $\Gamma, \iota$ and cut off at $T$ satisfies \eqref{def:decay} for $s$ and \eqref{def:levy}. 
	\end{lemma}
	\begin{proof}
		Fix $x_0 \in \Rd$ and $r > 0$. We observe that $\constM  < \constW |t|^\constExp$ if and only if $\big(\constM/\constW\big)^\frac{1}{\constExp} < |t|$. Then
		\begin{align*}
			\rho\big(x_0, B_r(x_0)^c\big)
			&= \int_{-T}^T \haus^0\big(\{y \in B_r(x_0)^c \mid \viewMap_t(y) = x_0\}\big) \, \iota(\intd t)\\
			&\leq c (\constI) \Big[
			\int_{r}^{(\constM/\constW)^{\frac{1}{\constExp}} \wedge r} \constM t^{-1-2s} \, \intd t
			+ \int_{(\constM/\constW)^{\frac{1}{\constExp}} \wedge r}^T
			\constW t^{-1-2s+\constExp} \, \intd t
			\Big]\\
			&< c(\constI) \Big[\frac{\constM}{2s} r^{-2s}
			+ \frac{\constW T^\constExp}{2s} r^{-2s}\Big],
		\end{align*}
		where we have used \eqref{ass:w} and \eqref{def:iota-comparable}. This proves \eqref{def:decay}.
		The verification of \eqref{def:levy} works in a similar way, see \cite[Lemma 4.3.7]{HeppDiss}.
	\end{proof}
	
	We can now prove \autoref{thm:properties-for-regularity-theory}.
	\begin{proof}[Proof of \autoref{thm:properties-for-regularity-theory}]
		
		Let $T_0$ be the maximal cutoff point from \autoref{measures:lemma:rearview-measure-is-well-defined} (b) such that the family of measures $\rearview(x,  \intd y)$ given in \autoref{meas:def:rearview-measure} is well-defined. In particular, $(x, t) \mapsto \haus^0\big(\{y \in \Rd \mid \viewMap_t(y) = x\}\big)$ will be $\leb \big(\Rd \times (-T_0, T_0)\big)$-measurable. 
		The family of measures $\fullview(x, \intd y)$ in \eqref{def:fullview-concrete} will then be constructed as the sum of $\frontview(x, \intd y)$ as given in \autoref{def:frontview-concrete} and $\rearview(x, \intd y)$, both with cutoff point $T \in (0, T_0)$. Well-definedness of $\fullview(x, \intd y)$ follows from \autoref{measures:lemma:nonsymmetric-measure-well-defined} and \autoref{measures:lemma:rearview-measure-is-well-defined}. 
		
		Let $A, B \in \borel(\Rd)$. We use the area formula (see e.g. \cite[Theorem 3.2.3]{Fed69}) and Tonelli's theorem to see that  
		\begin{align*}
			\int_A \frontview(x, B) \, \intd x
			&= \int_A \int_{-T}^T 
			\indic_B \big(x + t \Gamma(x)\big)
			|\det D\viewMap_t (x)| \, 
			\iota(\intd t) \intd x\\
			&= \int_{-T}^T \int_A 
			\indic_B \big(x + t \Gamma(x)\big)
			|\det D\viewMap_t (x)| \, 
			\intd x\, \iota(\intd t)\\
			&= \int_{-T}^T 
			\int_\Rd \indic_B(y) 
			\haus^0\big(\{ x \in A \mid \viewMap_t(x) = y\}\big)
			\, \intd y 
			\, \iota(\intd t)
			= \int_B \rearview(y, A) \, \intd y. 
		\end{align*}
		This implies \eqref{def:symmetry}, i.e. (a), even for intersecting sets. Note that in this case both integrals will be infinite.
		
		The above gives for a symmetric function $f(x, y) \geq 0$, 
		\begin{align*}
			\int_\domain \int_\domain f(x, y)\, \fullview(x, \intd y) \, \intd x 
			= \int_\domain \int_\domain f(x,y)\,\frontview(x, \intd y) \, \intd x,
		\end{align*}
		which, in turn, implies (c). 
		Under the additional assumption \eqref{ass:w}, (b) follows from \autoref{measures:lemma:nonsymmetric-measure-well-defined} (b) and \autoref{measures:lemma:rearview-decay-levy}. 
	\end{proof}

\section{Connectedness}\label{sec:connectedness}
	In this section, we show that it is possible to choose small cubes where all points are connected via at most $d$ jumps, provided that we are dealing with a Lipschitz continuous configuration. 
	\begin{theorem}[Connectedness on small scales]
	\label{comp:thm:connectedness-on-small-scales}
		Let $\constL> 0$ and $\constDet \in (0, 1)$. 
		There are a radius $R_0 = R_0(d, \constL, \constDet) \in (0, 1)$ and two factors $\lambda = \lambda(d, \constL, \constDet), \kappa=\kappa(d, \constL, \constDet)$ $> 1$, such that,
		for any configuration $\Gamma: \Rd \to (\Sd)^d$
		that satisfies \eqref{ass:lip} with Lipschitz constant $\constL$ and \eqref{ass:det} with constant $\constDet$,  
		as well as for any $R \in (0, R_0)$ and $z \in \Rd$, there  holds:
		\begin{enumerate}[(a)]
			\item For any $x \in B_R(z)$, $\pathMap_x: B_{\kappa R} \to \Rd$ is a bi-Lipschitz map onto its image with $B_R(z) \subset \pathMap_x(B_{\kappa R})$. Moreover, for any $y \in B_R(z)$, $\frac{2}{3} < |\det D \pathMap_x\inv(y) | < \frac{2}{\constDet}$ whenever $D \pathMap\inv_x(y)$ exists. 
			\item All points $x, y \in B_R(z)$ are connected. Moreover, it is possible to choose a path $\pathP(x, y)$ between $x$ and $y$ with the following properties: 
			\begin{enumerate}[(b.1)]
				\item The path has length at most $d$.
				\item The edges of the path are given by the path map $\pathMap_x$, i.e. 
				\begin{align*}
					\pathP(x, y)
					= 
					\big(\, p_k\big(x, \pathMap_x\inv(y)\big)\,\big)_{k=1}^d,
				\end{align*}
				where $p_k\big(x, \pathMap_x\inv(y)\big)$ is defined as in \eqref{comp:def:connecting-points}.
				It is also possible to give the connecting points in terms of $\pathMap_y$. 
				%\item All connecting points $\pathPoint_k(x, y)$ lie in a ball $B_{\lambda R}(z)$. 
				\item For any $x \in B_R(z)$, $\pathMap_x(B_{\kappa R}) \subset B_{\lambda R}(z)$. In particular, the connecting points all lie in the ball $B_{\lambda R}(z)$.
			\end{enumerate}
		\end{enumerate}
	\end{theorem}
	We start by introducing the path map $\pathMap_x$ for each point $x$, which assigns to a vector $q = (q_1, \dots, q_d)$ the point that may be reached starting from $x$ and successively jumping $q_k$-far in direction $\Gamma_k$. 
	The ingredients of the proof of \autoref{comp:thm:connectedness-on-small-scales} then come from the two main assumptions on the configuration. Thanks to \eqref{ass:lip}, the path map is invertible on small scales. And due to \eqref{ass:det}, the image of the path map contains a small ball around $x$. In other words, a small ball around any point may be reached in at most $d$ jumps, starting from the point itself.  
	
	Let $x \in \Rd$ be a point and $q = (q_1, \dots, q_d) \in \Rd$ a vector of jump lengths. 
	A path from $x$ to a point $p_d(x, q)$ is described recursively by the connecting points
	\begin{align}\label{comp:def:connecting-points}
		p_0(x, q) &\defeq x 
		\text{ and }
		p_k(x, q) = q_{k-1}(x, q) + q_k \Gamma_k \big(p_{k-1}(x, q)\big),
	\end{align}
	for $k \in \{1, \dots, d\}$. 
	Note that two consecutive points in this chain are directly connected, and we write
	$\pathP\big(x, p_d(x, q)\big) = \big(p_k(x, q)\big)_{k=0}^d$ for the path. 
	For $x \in \Rd$, define the \textit{path map}
	\begin{align}\label{comp:def:path-map}\tag{PM}
		\pathMap_x: \Rd &\to \Rd\\
		q &\mapsto x + q_1\Gamma_1(x) + \dots + q_d \Gamma_d\big(p_{d-1}(x, q)\big)= p_d(x, q).\notag
	\end{align}
	If we can find a radius $R > 0$ and factors $\lambda, \kappa > 1$ such that $\pathMap_x: B_{\kappa R} \to B_{\lambda \kappa R}(x)$ is invertible with $B_R(x) \subset \pathMap_x(B_{\kappa R})$, all points in $B_R(x)$ must be connected to $x$.
	The path will be given by $\pathP \big(x, \pathMap_x(q)\big)$ and all intermediate points in the path will lie in the larger ball $B_{\lambda \kappa R}(x)$. 
%	This connectedness result on small scales, together with the characterization of the paths, will allow us to deduce coercivity on small scales. 

	Observe that, for each $x \in \Rd$, the path map $\pathMap_x$ is locally Lipschitz on all of $\Rd$ (see \cite[Lemma 5.1.2]{HeppDiss}).
	Moreover, all connecting points between $x$ and $p_d(x, q)$ for some $q \in B_R$ lie in a slightly larger ball around $x$ whose radius is a multiple of $R$:
	\begin{lemma}\label{comp:lemma:connecting-points-lie-in-slightly-bigger-ball}
		There is a constant $\lambda = \lambda(d) > 1$ such that, for any configuration $\Gamma: \Rd \to (\Sd)^d$, $R > 0$, $x \in \Rd$, and $q \in B_R$, there holds $p_k(x, q) \in B_{\lambda R}(x)$ for any $k \in \{1, \dots, d\}$.
		In particular, $\pathMap_x(B_R) \subset B_{\lambda R}(x)$. 
	\end{lemma}
	\begin{proof}
		Observe that $\pathMap_x(0) = x$. Let $q \in B_R$ and $k \in \{1, \dots, d\}$. Then 
		\begin{align*}
			\big|x - p_k(x, q)\big| 
			&= \Big|x - \Big(x + \sum_{j=1}^k q_j \Gamma_j\big(p_{j-1}(x, q)\big)\Big)\Big|
			\leq \sum_{j=1}^k |q_j|\big|\Gamma_j\big(p_{j-1}(x, q)\big)\big|
			\leq d R, 
		\end{align*}
		since $|q_j| \leq |q| < R$ by assumption and $|\Gamma_j(z)| = 1$ for any $z \in \Rd$. 
	\end{proof}
	%It is possible to obtain a better constant if we assume all type vector matrices $\RotM_x$ to be very close to orthogonal matrices. However, once we allow type vectors belonging to the same point to be more and more similar to each other, the constant becomes increasingly worse.  
	
	By Rademacher's theorem, the partial derivatives of each subconfiguration $\Gamma_k$ exist almost everywhere and the same holds for each path map $\pathMap_x$.
	An important technical result is the following:
	\begin{proposition}\label{comp:prop:norm-of-jacobian-of-path-map}
		Let $\delta, \constL > 0$. 
		There is a radius $R = R(d, \constL, \delta) > 0$, such that,
		for any configuration $\Gamma = \Rd \to (\Sd)^d$ that satisfies 
		\eqref{ass:lip} with Lipschitz constant $\constL$, 
		as well as for any $x \in \Rd$ and $q \in B_R$ such that $D\pathMap_x(q)$ exists, there holds
		\begin{align*}
			|| D \pathMap_x (q) - \RotM_x||_F < \delta. 
		\end{align*}
	\end{proposition}
	The full proof can be found in \cite[Proposition 5.1.5]{HeppDiss}. We merely summarize the strategy here. Using the chain rule, one calculates the partial derivatives of $\pathMap_x$ to obtain
	\begin{align*}
		\big(D \pathMap_x(q)\big)_{ij}
		= \Gamma_j\upperInd{i}\big(p_k(x, q)\big)
		+ \sum_{k= j+1}^d q_k R(\partial_k \Gamma\upperInd{i}, q),
	\end{align*}
	where the remainder term depends on the values of the partial derivatives of $\Gamma$. For a $d$-tuple of points $y\upperInd{1}, \dots, y\upperInd{d}$ in $\Rd$, set 
	\begin{align}\label{comp:eq:notation-for-type-matrix-with-different-points}
		B_{\{y\upperInd{1}, \dots, y\upperInd{d}\}}
		\defeq 
		\begin{pmatrix}
			| 	&	& | \\
			\Gamma_1(y\upperInd{1}) 	& \cdots 	& \Gamma_d(y\upperInd{d})\\
			|	&	& |
		\end{pmatrix}
	\end{align}
	and observe that $\RotM_x$ is close to $B_{\{y\upperInd{1}, \dots, y\upperInd{d}\}}$ if all points $y\upperInd{k}$ are close to $x$, due to \eqref{ass:lip}. 
	Since all points $p_k(x, q)$ are close to $x$ if $q$ is small, thanks to \autoref{comp:lemma:connecting-points-lie-in-slightly-bigger-ball}, one only needs to estimate the remainder term by 
	\begin{align*}
		\Big|\sum_{k= j+1}^d q_k R(\partial_k \Gamma\upperInd{i}, q)\Big|
		\leq c R \constL
		\text{, if } |q| < R, 
	\end{align*}
	to obtain \autoref{comp:prop:norm-of-jacobian-of-path-map}. 
	
	In order to show invertibility of the path map on small balls, we prove a version of the inverse function theorem for Lipschitz functions, see \cite[Theorem 1]{Clarke1976} and \cite[Chapter 7]{ClarkeBook}. Under significantly stronger assumptions on the given Lipschitz map, it is possible to quantify the domain on which the function can be locally inverted. 
	The proof uses a special case of a result about the Neumann series of a bounded linear operator in a Banach space (see \cite[Theorem II.1.1]{Wer2018}).
	The result states that if $A \in \R^{d \times d}$ is a matrix with $||A||_{op} < 1$, then $(\id - A)$ is invertible and $||(\id - A)\inv||_{op} \leq (1 - ||A||_{op})\inv$. Recall that $\id$ denotes the identity matrix in $\R^{d\times d}$. 
	We deduce the following remark, proved in \cite[Remark 5.1.7]{HeppDiss}.
	\begin{remark}\label{comp:rem:Neumann-series}
		Let $S \in \R^{d\times d}$ be an invertible matrix and $A \in \R^{d\times d}$ be such that $||S - A ||_{op} < \delta ||S\inv||_{op}\inv$ for some $\delta \in (0, 1)$. Then $A$ is invertible as well and $||A\inv||_{op} \leq \frac{1}{1-\delta}||S\inv||_{op}$. 
	\end{remark}
%	\begin{proof}
%		Write $A = S \big(\id - (\id - S\inv A)\big)$. Then 
%		\begin{align*}
%			||\id - S\inv A||_{op} 
%			= ||S\inv(S - A)||_{op}
%			\leq ||S\inv||_{op}||S-A||_{op} < \delta < 1.
%		\end{align*}
%		Using the Neumann series, we see that $\id - (\id - S\inv A)$ is invertible with $||\big( \id - (\id - S\inv A)\big)\inv||_{op} \leq \big(1 - ||\id - S\inv A||_{op}\big)\inv$. Since $A$ is the product of invertible matrices, it is itself invertible and there holds
%		\begin{align*}
%			||A\inv||_{op} 
%			&= ||\big(\id - (\id-S\inv A)\big)\inv S\inv ||_{op}\\
%			&\leq \big(1 - ||\id - S\inv A||_{op}\big)\inv ||S\inv||_{op}  < \tfrac{1}{1-\delta}||S\inv||_{op}.
%		\end{align*}
%	\end{proof}
	
	The inverse function theorem for Lipschitz functions works with the concept of a generalized Jacobian, which will be introduced now. 
	The set of matrices $\R^{d \times d}$ shall have the topology of the Frobenius norm. 
	If $f: \Rd \to \Rd$ is Lipschitz around some point $x \in \Rd$, set 
	\begin{align*}
		J f(x) \defeq \big\{ \lim_{i\to\infty} D f(x_i) 
		\mid Df(x_i) \text{ exists, and } x_i \to x \big\} \subset \R^{d \times d}.
	\end{align*}
	The \textit{generalized Jacobian} of $f$ is defined as the convex hull of this set, i.e. 
	\begin{align}\label{def:generalized-jacobian}
		\partial f (x) 
		\defeq 
		\conv Jf(x) \subset \R^{d \times d}.
	\end{align}
	If every $M \in \partial f(x)$ is invertible, we say that $\partial f(x)$ is of \textit{maximal rank} (cf. \cite[Definitions 1 and 2]{Clarke1976}). 
	We derive a characterization of the generalized Jacobian of a Lipschitz function under quite strong conditions on the differentials in each point. 
	\begin{lemma}\label{lemma:generlized-jacobian-full-rank}
		Let $c_0 > 0$, $\delta \in (0, 1)$, and $f: \Rd \to \Rd$ be a function that is Lipschitz continuous on some open neighborhood $U \subset \Rd$. 
		Suppose there is some invertible matrix $S \in \R^{d \times d}$ with $||S\inv||_{op} < c_0$, such that every $x_0 \in U$ where $Df(x_0)$ exists satisfies 
		\begin{align}
			\label{comp:eq:inverse-function-theorem-ass}
			||Df(x_0)- S ||_{op} < \delta c_0\inv.
		\end{align}
		Then there holds 
		\begin{align*} 
			\forall \; x_0 \in U, M \in \partial f(x_0) \; : \; 
			M \text{ is invertible with }
			||M\inv||_{op} < \frac{c_0}{1 - \delta}. 
		\end{align*} 
		In particular, for each $x_0 \in U$, the generalized Jacobian $\partial f(x_0)$ is of full rank. 
	\end{lemma}
	\begin{proof}
		Let $x_0 \in U$. 
		If $A \in Jf(x_0)$, then $A = \lim_{i\to \infty} Df(x_i)$ for some $x_i \to x_0$. Since $U$ is open, we may assume \eqref{comp:eq:inverse-function-theorem-ass} to hold for each $Df(x_i)$. This gives 
		\begin{align*}
			||A - S||_{op} 
			= ||\lim_{i \to \infty} Df(x_i) - S||_{op} 
			= \lim_{i\to \infty} ||Df(x_i) - S||_{op}
			< \delta c_0\inv. 
		\end{align*}
		
		Now let $M \in \partial f(x_0)$, i.e. $M = \sum_{i=1}^n s_i A_i$, where $A_i \in Jf(x_0)$, $s_i \geq 0$, and $\sum_{i=1}^n s_i = 1$. Then 
		\begin{align*}
			||M - S||_{op} 
			= \Big|\Big|
			\sum_{i=1}^n s_i A_i - \sum_{i=1}^n s_i S\, 
			\Big|\Big|_{op}
			\leq \sum_{i=1}^n s_i ||A_i - S||_{op}
			< \delta c_0\inv. 
		\end{align*}
		Thus, $||M - S||_{op} < c_0$ holds for any matrix of the generalized Jacobian of $f$ at $x_0$. 
		\autoref{comp:rem:Neumann-series} together with $||S\inv||_{op} < c_0$ implies 
		\begin{align*}
			\forall \; M \in \partial f(x_0) \; : \; 
			M \text{ is invertible and }
			||M\inv||_{op} \leq \frac{1}{1-\delta} ||S\inv||_{op}
			< \frac{c_0}{1-\delta}.
		\end{align*}
	\end{proof} 
	
	\begin{proposition}[Quantified inverse function theorem for Lipschitz functions]\label{comp:prop:quantified-inverse-function-theorem}
		Let $c_0, \rho > 0$ and $\delta \in (0, 1)$. 
		Let $f: \Rd \to \Rd$ be a function that is Lipschitz continuous on some ball $B_\rho(x_0)$ for $x_0 \in \Rd$. 
		Assume furthermore that there is an invertible matrix $S \in \R^{d \times d}$ with $||S\inv||_{op} < c_0$, such that
		any $x \in B_{\rho}(x_0)$ where $D f(x)$ exists satisfies \eqref{comp:eq:inverse-function-theorem-ass}.
		Then there are a radius $r_0 = r_0(c_0, \delta) \in (0, \rho)$ and a factor $\kappa = \kappa (c_0, \delta) > 0$, such that any $r \in (0, r_0)$ satisfies
		\begin{enumerate}[(a)]
			\item $f: B_r(x_0) \to \Rd$ is injective, 
			\item $f\big(B_r(x_0)\big) \supset B_{\kappa r}\big(f(x_0)\big)$, 
			\item the inverse function $f\inv: B_{\kappa r}\big(f(x_0)\big) \to B_r(x_0)$ is Lip\-schitz continuous with Lip\-schitz constant depending only on $\delta$ and $c_0$. 
		\end{enumerate}
	\end{proposition}
	\begin{proof}
		The proof of \autoref{comp:prop:quantified-inverse-function-theorem} follows the proof of the inverse function theorem for Lipschitz functions in \cite[Theorem 1]{Clarke1976}, which is split into several lemmas. The inverse function theorem as well as all lemmas making up its proof are shown under the assumption that $f$ is Lipschitz around $x_0$ and of maximal rank. 
		By \autoref{lemma:generlized-jacobian-full-rank}, these assumptions are also true for any $f$ satisfying the assumptions of \autoref{comp:prop:quantified-inverse-function-theorem} at $x_0$. Thus, all lemmas in Clarke's proof of the inverse function theorem also hold for such a function $f$. 
		
		\autoref{comp:prop:quantified-inverse-function-theorem} even implies that any $M \in \partial f(x_0)$ is invertible with $||M\inv||_{op}< \frac{c_0}{1-\delta} =: c_1$. 
		This gives a stronger version of \cite[Lemma 3]{Clarke1976}: There are a radius $r = r(c_0, \delta) > 0$ and a factor $\kappa = \kappa(c_0, \delta) > 0$, such that, for $v \in \Sd$, there is a unit vector $w \in \Sd$, with
		\begin{align}\label{comp:eq:inverse-function-theorem-lemma3-claim}
			\forall \; x \in B_r(x_0), M \in \partial f(x_0) 
			\; : \; 
			w \cdot (Mv) \geq \kappa. 
		\end{align}
		The main difference to \cite[Lemma 3]{Clarke1976} are the dependencies of $r$ and $\kappa$. Due to the weaker assumptions, Clarke can only conclude the existence of \textit{some} $r$ and $\kappa$, with their value depending on $x_0$. In our case, the values of $r$ and $\kappa$ will be the same for any $x_0$ satisfying the theorem's assumptions, since they can be chosen only depending on $\delta$ and $c_0$. 
		
		In order to prove the lemma, we consider $\partial f(x_0) \Sd \defeq \{Mv \in \Rd \mid M \in \partial f(x_0), v \in \Sd\}$. 
		Observe that $\dist \big(0, \partial f(x_0)\Sd\big) > \frac{1}{c_1}$. 
		Indeed, if $y = Mv \in \partial f(x_0)\Sd$, then 
		\begin{align*}
			|y| 
			= |Mv| 
			\geq \frac{|v|}{||M\inv||_{op}}
			> \frac{1}{c_1}. 
		\end{align*}
		Next, let $G \defeq \partial f(x_0) + B^{d \times d}_{\frac{1}{2c_1}} \subset \R^{d \times d}$, where $B^{d \times d}_{\frac{1}{2c_1}}$ is the ball of radius $\frac{1}{2c_1}$ centered at $0$ in $\R^{d\times d}$. Again, we can bound the distance between $G\Sd$ and $0$ from below, obtaining  $\dist(0, G\Sd) > \frac{1}{2 c_1}$. Indeed, $|Av| \leq ||A||_{op} \leq ||A||_F < \frac{1}{2c_1}$ if $A \in B^{d \times d}_{\frac{1}{2c_1}}$. For $y \in G\Sd$, we have $y = Mv + Av$ with $M \in \partial f(x_0)$, $A \in B^{d \times d}_{\frac{1}{2c_1}}$, which implies 
		\begin{align*}
			|Mv + Av| 
			\geq |Mv| - |Av| 
			> \frac{1}{c_1} - \frac{1}{2c_1}
			= \frac{1}{2 c_1}. 
		\end{align*}
		
		According to \cite[Lemma 1]{Clarke1976}, for each $\eps > 0$, there is $r_0 = r_0(\eps) > 0$, such that $\partial f(x) \subset \partial f(x_0) + \eps B_1^{d \times d}$, if $x \in B_{r_0}(x_0)$. Naturally, we also need to assume $r_0 < \rho$ so that $f$ is still Lipschitz on $B_{r_0}(x_0)$. 
		In particular, if $\eps = \frac{1}{2 c_1}$, the radius $r_0$ will only depend on $c_1$, i.e. $c_0$ and $\delta$, and thus $\partial f(x) \subset G$ if $x \in B_{r_0}(x_0)$. 
		
		If $v \in \Sd$, $Gv$ is convex and, by the above, $\dist(0, Gv) > \frac{1}{2c_1}$. By the hyperplane separation theorem, there is then some $w \in \Sd$ such that any $M \in G$ verifies $w \cdot (Mv) \geq \frac{1}{2c_1}$. The conclusion follows from $\partial f(x) \subset G$.
		Now $r_0$ is the desired radius and $\kappa = \frac{1}{2c_1}$. 
		
		By \cite[Lemma 4]{Clarke1976}, there holds 
		\begin{align*}
			\forall \; x, y \in \overline{B_{r_0}(x_0)}
			\; : \; 
			|f(x) - f(y)| \geq \kappa |x-y|.
		\end{align*}
		Hence $f$ is injective on $B_r(x_0)$. Moreover, due to \cite[Lemma 5]{Clarke1976}, $f\big(B_{r_0}(x_0)\big) \supset B_{\kappa r_0}\big(f(x_0)\big)$. 
		The inverse function $g$ of $f$ is then defined on $B_{\kappa r_0}\big(f(x_0)\big)$ as $g(y) = x$, where $x$ is the unique point in $B_{r_0}(x_0)$ such that $f(x) = y$. Due to the injectivity of $f$, this assignment is well-defined. Moreover, $g$ is Lipschitz with constant $\kappa$, due to the above inequality for $f$. 
		The same considerations apply to any radius $r \in (0, r_0)$. 
	\end{proof}

	Let us now prove connectedness in a small neighborhood of a point. 
	\begin{corollary}\label{comp:cor:connectedness-around-point}
		Let $\constL> 0$ and $\constDet \in (0, 1)$. 
		There are a radius $R_0 = R_0(d, \constL, \constDet) \in (0, 1)$ and two factors $\lambda = \lambda(d), \kappa=\kappa(d, \constL, \constDet) > 1$, such that,
		for any configuration $\Gamma: \Rd \to (\Sd)^d$
		that satisfies \eqref{ass:lip} with Lipschitz constant $\constL$ and \eqref{ass:det} with constant $\constDet$,  
		as well as for any $R \in (0, R_0)$ and $x \in \Rd$, there  holds:
		\begin{enumerate}[(a)]
			\item The path map $\pathMap_x: B_{\kappa R} \to \Rd$ is bi-Lipschitz onto its image. 
			Moreover, $B_R(x) \subset \pathMap_x(B_{\kappa R})$ and $\frac{2}{3} < |\det D \pathMap_x\inv(y)| < \frac{2}{\constDet}$ for every $y \in B_R(x)$ for which $D \pathMap_x\inv(y)$ exists. 
			\item All points in $B_R(x)$ are connected with $x$ via a path that has $d+1$ edges. 
			Specifically, if $y \in B_R(x)$, then the path between $x$ and $y$ is given by
			\begin{align*}
				\pathP (x, y) 
				= \big(\, p_k\big(x, \pathMap_x\inv(y)\big)\, \big)_{k=0}^d \text{ with } p_k\big(x, \pathMap_x\inv(y)\big) \text{ as in \eqref{comp:def:connecting-points}.}
			\end{align*}
			\item The edges $p_k\big(x, \pathMap_x\inv(y)\big)$ of the path $\pathP(x, y)$ all lie in $B_{\lambda \kappa R}(x)$. 
		\end{enumerate}
	\end{corollary}
	\begin{proof}
		Due to \autoref{comp:prop:norm-of-jacobian-of-path-map} paired with the continuity of the determinant, there is a radius $R_1 = R_1(d, \constL, \constDet) \in (0, 1)$, such that, for any $q \in B_{R_1}$ where $D \pathMap_x(q)$ exists, there holds $\frac{\constDet}{2} < |\det D \pathMap_x(q)| < \frac{3}{2}$. 
		Next, we verify the assumptions of the quantified inverse function theorem  \autoref{comp:prop:quantified-inverse-function-theorem}. 
		We know that $\pathMap_x$ is Lipschitz on $B_{R_1}$. 		
		By \eqref{meas:cc:config-det:eq:norm-of-inverse-bound}
		the type vector matrix $\RotM_x$ at $x$ is invertible and there even is a constant $c_0 = c_0(d, \constDet) > 0$ such that $||\RotM_x\inv||_{op} < c_0$.
		Let $\delta = \frac{1}{2}$. By \autoref{comp:prop:norm-of-jacobian-of-path-map}, there is a radius $R_2 = R_2(d, \constL, \constDet) \leq R_1$, such that, again for any $q \in B_{R_2}$ where $D \pathMap_x(q)$ exists, $||D \pathMap_x(q) - \RotM_x||_{op} < \frac{1}{2} c_0\inv$. 
		Thus all requirements for  \autoref{comp:prop:quantified-inverse-function-theorem} are met.
		This gives us a radius $R_3 = R_3(d, \constL, \constDet) < R_2$ and a constant $\tilde \kappa = \tilde \kappa (d, \constL, \constDet) > 0$ such that $\pathMap_x$ is injective on $B_{R_3}$ with $B_{\tilde \kappa R_3}\big(\pathMap_x(0)\big)  \subset \pathMap_x(B_{R_3})$. 
		Since $c_0 > 1$, $\tilde \kappa < 1$. Now set $\kappa \defeq \frac{1}{\tilde \kappa}$ and $R_0 \defeq \frac{R_3}{\kappa}$. Injectivity of $\pathMap_x$ on $B_{\kappa R_0}$ and the containment of $B_{R_0}(x) = B_R\big(\pathMap_x(0)\big)$ in $\pathMap_x(B_{\kappa R_0})$ follow immediately. In particular, $D \pathMap_x\inv(y)$ exists for almost every $y \in B_{R_0}(x)$ and the bounds for $|\det D \pathMap_x\inv|$ follow from the bounds for $|\det D \pathMap_x|$. This proves (a).
		
		Further, by (a), if $y \in B_R(x)$, there is some $q \in B_{\kappa R}$ with $\pathMap_x(q) = y$. Thus $\big(p_k(x, q) = p_k\big(x, \pathMap_x\inv(y)\big)\big)_{k=0}^d$ forms a path between $x$ and $y$, and (b) follows. 		
		Finally, by \autoref{comp:lemma:connecting-points-lie-in-slightly-bigger-ball}, each $p_k(x, q)$ lies in $B_{\lambda \kappa R}(x)$ for any $q \in B_{\kappa R}$, which gives (c).
	\end{proof}
	
	\begin{proof}[Proof of \autoref{comp:thm:connectedness-on-small-scales}]
		Let $R_1 > 0$ be the radius and $\tilde \kappa = \tilde \kappa(d, \constL, \constDet)$ the factor from \autoref{comp:cor:connectedness-around-point} and set $R_0 \defeq \frac{R_1}{2}$. If $R \in (0, R_0)$ and $x \in B_R(z)$, then $B_{2R}(x) \supset B_R(z)$. Thus $\pathMap_x$ is bi-Lipschitz on $B_{\tilde \kappa 2R}$ with $B_R(z) \subset B_{2R}(x) \subset \pathMap_x (B_{\tilde \kappa 2r})$. Set $\kappa \defeq 2 \tilde \kappa$. 
		Now \autoref{comp:cor:connectedness-around-point} gives (a), as well as (b.1) and (b.2). If $\tilde \lambda = \tilde \lambda (d)$ is the other factor from \autoref{comp:cor:connectedness-around-point}, then (b.3) follows for $\lambda \defeq (1 + \kappa \tilde \lambda)$ because
		$\pathMap_x (B_{\kappa R})
		\subset B_{\tilde \lambda \kappa R}(x)
		\subset B_{\tilde \lambda \kappa R + R}(z)$. 			
	\end{proof}
	
	The difference between \autoref{comp:thm:connectedness-on-small-scales} and \autoref{comp:cor:connectedness-around-point} might escape a cursory glance. In \autoref{comp:cor:connectedness-around-point}, it is proved that all points in a ball are connected to its center, while \autoref{comp:thm:connectedness-on-small-scales} shows that two arbitrary points in some (small) ball are connected \textit{to each other}. This is the result we will be using to prove the lower bound. 
		
	From \autoref{comp:thm:connectedness-on-small-scales}, we immediately gain a connectedness result on the whole space by covering $\Rd$ with balls of a small but fixed radius. 
	\begin{corollary}
		Let $\Gamma: \Rd \to (\Sd)^d$ be a configuration that satisfies \eqref{ass:det} and \eqref{ass:lip}. Then all points in $\Rd$ are connected. 
	\end{corollary}

\section{Coercivity}\label{sec:coerc}
	In this section, we prove the lower bound of \eqref{eq:comparability}:
	\begin{theorem}[Coercivity for Lipschitz configurations]\label{thm:coercivity-lip}
		Let $s, \constDet \in (0, 1)$, $\constI> 1$, and $\constL, T_0> 0$. If $\domain \subset \Rd$ is an open, bounded, and connected Lipschitz domain, 
		there is a constant $c = c(d, s, \domain, \constI, \constL, \constDet, T_0) > 0$, such that, 
		for any measurable function $u: \domain \to \R$, 
		\begin{align}
			\quadForm^s_\domain (u, u)
			\, \leq \, c \, \quadForm^\Gamma_\domain(u, u), 
		\end{align}
		where $\quadForm^\Gamma_\domain$ is any quadratic form with cutoff point $T \in (T_0, \infty]$ generated by 
		a configuration $\Gamma: \Rd \to (\Sd)^d$ that satisfies
		\eqref{ass:lip} with Lipschitz constant $\constL$ 
		and 
		\eqref{ass:det} with constant $\constDet$,
		and by a one-dimensional Lebesgue measure $\iota(\intd t)$ that satisfies \eqref{def:iota-comparable} with $\constI$ for $s$.
	\end{theorem}
	%The main ingredient of this theorem's proof is the connectedness result \autoref{comp:thm:connectedness-on-small-scales} from the previous section. 
	
\subsection{Coercivity on small scales}\label{sec:coerc:small-scales}
	Let $q \in \Rd$ and $k \in \{1, \dots, d\}$. The \textit{dual path map} is given by 
	\begin{align}%\label{comp:def:dual-path-map}
		\dPathMap\upperInd{k}_q : \Rd &\to \Rd \\%\tag{DPM}\\
		x &\mapsto p_k(x, q).\notag
	\end{align}
	In this setup, a vector $q$ of jump lengths is held fixed and the starting point $x$ is mapped to the $k$-th point in its path to $p_d(x, q) = \pathMap_x(q)$. 
	We have the useful characterization 
	\begin{align}\label{comp:eq:dual-path-map-characterization} 
		\dPathMap_q\upperInd{k+1}(x) 
		= \dPathMap_q\upperInd{k}(x) + q_{k+1} \Gamma_{k+1} \big( \dPathMap_q\upperInd{k}(x)\big)
		= p_k(x, q) + q_{k+1} \Gamma_{k+1} \big(p_k(x, q)\big). 
	\end{align}
	Lipschitz continuity of each $\dPathMap_q\upperInd{k}$ follows from lemma \cite[Lemma 5.2.2]{HeppDiss}. 
	In the same way as before, it is possible to show 	
	\begin{align*}
		||D \dPathMap_q\upperInd{k}(x) - \id ||_F < \delta.
	\end{align*}
	for any $x$, granted that $q$ is small enough, see \cite[Lemma 5.2.4]{HeppDiss}.
%	We now prove an analog of \autoref{comp:prop:norm-of-jacobian-of-path-map} for $\dPathMap_q\upperInd{k}$. 
%	\begin{lemma}\label{comp:lemma:dual-path-map-norm-jacobian}
%		Let $\constL, \delta > 0$. 
%		There is a radius $R = R(d, \constL, \delta) \in (0, 1)$, 
%		such that, for any configuration $\Gamma: \Rd \to (\Sd)^d$ 
%		that satisfies \eqref{ass:lip} with Lipschitz constant $\constL$, 
%		as well as for any $q \in B_R$, $k \in \{1, \dots, d\}$, and every $x\in \Rd$ for which $D \dPathMap_q\upperInd{k}(x)$ exists, there holds 
%		\begin{align*}
%			||D \dPathMap_q\upperInd{k}(x) - \id ||_F < \delta.
%		\end{align*}
%	\end{lemma}
	
	From this, we deduce, using the continuity of the determinant as well as the quantified inverse function theorem (\autoref{comp:prop:quantified-inverse-function-theorem}):
	\begin{proposition}\label{comp:prop:dual-path-map-diffeo}
		Let $\constL, \kappa > 0$ and $\eps \in (0, 1)$. 
		There are a radius $R_0 = R_0(d, \constL, \kappa , \eps) \in (0, 1)$ and a factor $\lambda = \lambda(d, \kappa) > 1$, such that, 
		for any configuration $\Gamma: \Rd \to (\Sd)^d$
		that satisfies \eqref{ass:lip} with Lipschitz constant $\constL $, 
		as well as for any $R \in (0, R_0)$, $q \in B_{\kappa R}$, $k \in \{1, \dots, d\}$, and $z \in \Rd$, there holds
		\begin{enumerate}[(a)]
			\item for any $x \in \Rd$, if $D \dPathMap_q\upperInd{k}(x)$ exists, then $1- \eps < |\det D \dPathMap_q\upperInd{k}(x)| < 1 + \eps$,
			\item $\dPathMap_q\upperInd{k}: B_{R}(z) \to \Rd$ is a bi-Lipschitz map onto its image, and
			\item $\dPathMap_q\upperInd{k}\big(B_R(z)\big) \subset B_{\lambda R}(z)$. 
		\end{enumerate} 
	\end{proposition}
	See \cite[Proposition 5.2.5]{HeppDiss} for a proof. 
	Before we show coercivity on small scales, we prove a technical lemma that expresses the intuition that $|q|$ is comparable to the distance between $x$ and its image under $\pathMap_x$, if it is small enough. 
	\begin{lemma}\label{comp:lemma:comparability-norm-estimate}
		Let $\constL> 0$, $\constDet \in (0, 1)$. 
		There are a radius $R = R(d, \constL, \constDet) > 0$ and a 
		constant $c = c(d, \constDet) > 0$, 
		such that, for any configuration $\Gamma: \Rd \to (\Sd)^d$ that satisfies 
		\eqref{ass:lip} with Lipschitz constant $\constL$ and 
		\eqref{ass:det} with constant $\constDet$, 
		as well as for any $x \in \Rd$, and $q \in B_R$, there holds 
		$	|q| 
			\leq c
			\big| \pathMap_x(q) - x\big|.$
	\end{lemma}
	
	\begin{proof}
		Let $\lambda = \lambda(d)$ be the constant from \autoref{comp:lemma:connecting-points-lie-in-slightly-bigger-ball} and $c_0 = c_0(d, \constDet)>0$ the constant from \eqref{meas:cc:config-det:eq:norm-of-inverse-bound}, such that $||\RotM_x\inv||_{op} < c_0$. 
		For $q \in \Rd$, let $B_q \defeq B_{\{y\upperInd{1}, \dots, y\upperInd{d}\}}$ with $y\upperInd{k} = p_k(x, q)$ in the notation \eqref{comp:eq:notation-for-type-matrix-with-different-points}. 
		If $q \in B_R$, then $p_k(x, q) \in B_{\lambda R}(x)$ due to \autoref{comp:lemma:connecting-points-lie-in-slightly-bigger-ball}. 
		Thus, we may choose $R = R(d, \constL, \constDet)$ small enough such that $||B_q - \RotM_x||_{op} < \frac{1}{2}c_0\inv < \frac{1}{2}||\RotM_x\inv||_{op}\inv$. 
		Then $B_q$ is invertible with $||B_q\inv||_{op} < 2 ||\RotM_x\inv||_{op} < 2 c_0$, due to  \autoref{comp:rem:Neumann-series}. 
		
		So fix $q \in B_R$. 
		Observe that $\sum_{k=1}^d q_k \Gamma_k \big(p_{k-1}(x, q)\big) \in \Rd$ is a vector written in terms of the basis $\Gamma_1\big(p_1(x, q)\big), \dots, \Gamma_d\big(p_d(x, q)\big)$, whose base change matrix is precisely $B_q\inv$. 
		Since $\pathMap_x(q) - x 
		= \sum_{k=1}^d q_k \Gamma_k \big(p_{k-1}(x, q)\big)$, $B_q\inv\big(\pathMap_x(q) - x\big) = q$ and therefore 
		\begin{align*}
			|q| 
			\leq ||B_q\inv||_{op} \big|\pathMap_x(q) - x\big| 
			\leq 2 c_0 \big|\pathMap_x(q) - x\big|. 
		\end{align*}
	\end{proof}

	\begin{theorem}\label{comp:thm:lower-bound-for-singular-measures-small-scales} 
		Let $s, \constDet \in (0, 1)$, $\constI > 1$, and $\constL, T_0 > 0$.
		There are a radius $R_0 = R_0(d, \constL, \constDet, T_0) > 0$, a factor $\lambda = \lambda (d, \constL, \constDet) > 1$, and a constant $c = c(d, s, \constI, \constL, \constDet) > 0$, such that, for any $B_{\lambda R}(z) \subset \Rd$ with $R \in (0, R_0)$, and $u: B_R(z) \to \R$ measurable, there holds
		\begin{align*}
			\int_{B_R(z)} \int_{B_R(z)} &\uxy\\
			&\leq c
			\int_{B_{\lambda R}(z)}\int_{B_{\lambda R}(z)} \big( u(x) - u(y)\big)^2 \frontview(x, \intd y) \,  \intd x, 
		\end{align*}
		where $\frontview(x, \intd y)$ is any family of measures with cutoff point $T \in (T_0, \infty]$
		generated by a configuration $\Gamma: \Rd \to (\Sd)^d$ that satisfies
		\eqref{ass:lip} with Lipschitz constant $\constL$, and
		\eqref{ass:det} with constant $\constDet$,
		and by a one-dimensional Lebesgue measure $\iota(\intd t)$ that satisfies \eqref{def:iota-comparable} with $\constI$ for $s$.
	\end{theorem}
	
	\begin{proof}
		In order to use the connectedness result and the properties of the dual path map, we need to choose $R_0$ appropriately. 
		Let $R_1$ be the radius from \autoref{comp:thm:connectedness-on-small-scales} with factor $ \kappa_1$.  
		Let $R_2$ be the radius and $\lambda_2$ the factor from \autoref{comp:prop:dual-path-map-diffeo} for $\eps = \frac 12$ and $\kappa = \kappa_1$. 
		Let $R_3$ be the radius from \autoref{comp:lemma:comparability-norm-estimate}. 
		Set $\lambda \defeq \kappa_1 + \lambda_2$.
		Finally, let $R_4$ be the radius from \autoref{measures:lemma:properties-of-det-view-map} (c) for $\kappa = \frac 12$.
		Set $R_0 \defeq \min\{R_1, R_2, R_3, \frac{R_4}{\lambda}, \frac{T_0}{2\lambda}\}$. 
		%Observe that $R_0$ and $\lambda$ depend on $d, \constL$, and $\constDet$ as was required. 
		We may now apply \autoref{comp:thm:connectedness-on-small-scales}, \autoref{comp:prop:dual-path-map-diffeo}, and \autoref{comp:lemma:comparability-norm-estimate} to any ball $B_R(z)$, and \autoref{measures:lemma:properties-of-det-view-map} (c) to $B_{\lambda R}(z)$ if $R \in (0, R_0)$.
		
		First we make use of the connectedness result \autoref{comp:thm:connectedness-on-small-scales}. 
		If $x \in B_R(z)$, $\pathMap_x$ is invertible on $B_R(z)$ with $\pathMap_x\inv\big(B_R(z)\big) \subset B_{\kappa_1 R}$. 
		We may thus write $y = p_d\big(x, \pathMap_x\inv(y)\big)$ for any $y \in B_R(z)$ and perform a change of variables $q \defeq \pathMap_x\inv(y)$, where we use that $|\det D \pathMap_x\inv(y)|$ is bounded from below. Subsequently applying  \autoref{comp:lemma:comparability-norm-estimate} to the norm $|y-x| = |\pathMap_x(q) - x|$, we obtain 
		\begin{align*}
			&\int_{B_R(z)} \int_{B_R(z)} U(x, y) |y-x|^{-d-2s} \, \intd y \intd x \\
			\leq &\; c\int_{B_R(z)} \int_{B_R(z)} U\big(x, p_d(x, \pathMap_x\inv(y))\big) \; 
			\big|\pathMap_x \big(\pathMap_x\inv(y)\big)-x\big|^{-d-2s} 
			|\det D \pathMap_x\inv(y)| \, 
			\intd y \intd x\\
			= &\; c\; \int_{B_R(z)} \int_{\pathMap_x\inv(B_R(z))} U\big(x,  p_d(x, q)\big) \; |\pathMap_x(q) -x|^{-d-2s} \;\intd q \intd x\\
			\leq &\; c\; \int_{B_R(z)} \int_{\pathMap_x\inv(B_R(z))} U\big(x,  p_d(x, q)\big) \; |q|^{-d-2s} \;\intd q \intd x.
		\end{align*}
		Inserting the intermediate path points between $x$ and $p_d(x, q)$ and using $\pathMap_x\inv\big(B_R(z)\big) \subset B_{\kappa_1 R}$, we arrive at the estimate
		\begin{align*}
			&\int_{B_R(z)} \int_{B_R(z)} U(x, y) |y-x|^{-d-2s} \, \intd y \intd x \\
			\leq &\; c \; \sum_{k=1}^d \int_{B_R(z)} \int_{B_{\kappa_1 R}} U\big(p_{k-1}(x, q), p_k(x, q)\big) \; |q|^{-d-2s} \; \intd q \intd x.
		\end{align*} 
		The constant $c$ depends on $d, s, \constL$, and $\constDet$.
		
		Fix one $k \in \{1, \dots, d\}$. We use the characterization \eqref{comp:eq:dual-path-map-characterization} to rewrite $p_{k-1}(x, q)$ and $p_k(x, q)$ in terms of $\dPathMap_q\upperInd{k-1}(x)$. 
		After applying Tonelli's theorem, we change variables to $w \defeq \dPathMap_q\upperInd{k-1}(x)$ using \autoref{comp:prop:dual-path-map-diffeo}:
		\begin{align*}
			&\int_{B_R(z)} \int_{B_{\kappa_1 R}} U\big(p_{k-1}(x, q), p_k(x, q)\big) \; |q|^{-d-2s} \; \intd q \intd x\\
			= & \int_{B_{\kappa_1 R}} \int_{B_R(z)} U\big(\dPathMap_q\upperInd{k-1}(x), \dPathMap_q\upperInd{k-1}(x) +  q_k \Gamma_k\big(\dPathMap_q\upperInd{k-1}(x)\big)\big) \; |q|^{-d-2s} \; \intd x \intd q\\
			\leq & \; c \; \int_{B_{\kappa_1 R}} \int_{\dPathMap_q\upperInd{k-1}(B_R(z))} U(w, w +  q_k \Gamma_k(w)\big) \; |q|^{-d-2s} \; \intd w \intd q
		\end{align*}
		Due to \autoref{comp:prop:dual-path-map-diffeo} (c), $\dPathMap_q\upperInd{k-1}\big(B_R(z)\big) \subset B_{\lambda_2 R}(z)$. So we switch integrals again and then observe $B_{\kappa_1 R} \subset \big(-\kappa_1 R, \kappa_1 R\big)^d$. This gives: 
		\begin{align*}
			&\int_{B_R(z)} \int_{B_{\kappa_1 R}} U\big(p_{k-1}(x, q), p_k(x, q)\big) \; |q|^{-d-2s} \; \intd q \intd x\\
			\leq & \; \int_{B_{\lambda_2 R}(z)} \int_{B_{\kappa_1 R}} U(w, w +  q_k \Gamma_k(w)\big) \; |q|^{-d-2s} \; \intd q \intd w\\
			\leq& \;c \;  \int_{B_{\lambda_2R}(z)} \int_{-\kappa_1 R}^{\kappa_1 R} U\big(w, w +  q_k \Gamma_k(w)\big) \; \int_{\R^{d-1}}|q|^{-d-2s} \; \intd q^{(k)} \intd q_k \intd w\\
			\leq& \;c \;  \int_{B_{\lambda_2R}(z)} \int_{-\kappa_1 R}^{\kappa_1 R} U\big(w, w +  q_k \Gamma_k(w)\big) \; |q_k|^{-1-2s} \; \intd q_k \intd w
		\end{align*}
		In the third line, we write $q^{(k)} = (q_1, \dots, q_{k-1}, q_{k+1}, \dots, q_d)$, using Tonelli's theorem to rearrange the order of integration. 
		For the last inequality, we single out $q_k$ by estimating the integral over $|q|$ with respect to $q \upperInd{k}$, using a trick from \cite[proof of Theorem 6.2]{DyKa20}, see also 
		\cite[Lemma A.3]{HeppDiss}. 
		
		If $x \in B_{\lambda_2 R}(z)$ and $|t| < \kappa_1R$, $x + t \Gamma_k(x) \in B_{\lambda R}(z)$ by definition of $\lambda$. Now we use the fact that $\lambda R < R_4$ to artificially insert the determinant term with \autoref{measures:lemma:properties-of-det-view-map} (c): 
		\begin{align*}
			&\int_{B_{\lambda_2R}(z)} \int_{-\kappa_1 R}^{\kappa_1 R} U\big(x, x +  t \Gamma_k(x)\big) \; |t|^{-1-2s} \; \intd t \intd x\\
			\leq & \; c \int_{B_\lambda R(z)} \int_\R 
			U(x, x + t \Gamma(x)) 
			\indic_{B_{\lambda R}(z)} \big( x + t \Gamma_k(x)\big) 
			\, |\det D\viewMapKT(x)|
			|t|^{-1-2s}
			\, \intd t \intd x
		\end{align*}
		
		Since $2\lambda R < T_0 \leq T$, we obtain overall
		\begin{align*}
			&\int_{B_R(z)} \int_{B_R(z)} \uxy \\
			\leq & \; c \; \sum_{k=1}^d \int_{B_{\lambda R}(z)} \int_\R 
			U \big(x, x + t \Gamma_k(x)\big) 
			%\big(u(x) - u(x + t \Gamma_k(x))\big)^2
			%&\phantom{mmmmmmmmmmmmmmm}
			\indic_{B_{\lambda R}(z)} \big( x + t \Gamma_k(x)\big) 
			\, |\det D\viewMapKT(x)|
			|t|^{-1-2s}
			\, \intd t \intd x\\
			\leq & \; c \; \int_{B_{\lambda R}(z)} \int_{B_{\lambda R}(z)} \big(u(x)- u(y)\big)^2 \configMeas(x, \intd y) \intd x.
		\end{align*}
		The last inequality follows from \eqref{def:iota-comparable}. This proves the theorem. 
	\end{proof}

\subsection{Coercivity on Lipschitz domains}\label{sec:coerc:lipschitz-domains}
	In this subsection, we prove \autoref{thm:coercivity-lip}. 
	The theorem follows from an abstract result: A lower bound on small scales for the quadratic form $\quadForm^m$ associated with  any measure $m\big(\intd (x, y)\big)$ in the product space implies a lower bound for $\quadForm^m$ on arbitrary balls with a controlled growth factor on the right-hand side. So let us assume: 
	\begin{assumption}\label{comp:ass:comparability-on-small-scales}
		For the measure $m: \borel(\Rd \times \Rd) \to [0, \infty]$, there is a radius $\radAss >0$, a factor $\lambdaAss > 1$, and a constant $\constAss  > 0$, such that, for any cube $Q = Q_R(z)$ with $R \in (0, \radAss)$, $z \in \Rd$, and any measurable function $u: \Rd \to \R$, there holds
		\begin{align*}
			\quadForm_Q^s(u, u)
			&\leq \constAss \int_{\lambdaAss Q} \int_{\lambdaAss Q}
			\big(u(x) - u(y)\big)^2 m(\intd(x, y)).
		\end{align*} 
	\end{assumption}
	Then it follows: 
	\begin{lemma}\label{comp:lemma:comparability-in-balls}
		Let $m: \borel(\Rd \times \Rd) \to [0, \infty]$ be a measure that satisfies \autoref{comp:ass:comparability-on-small-scales} for cubes of radius $\radAss$, factor $\lambdaAss$, and constant $\constAss$. 
		For any growth factor $\Lambda > 1$, and radius $\rho > 0$, 
		there is a constant $c = c(d, s, \constAss, \radAss, \lambdaAss, \Lambda, \rho) > 0$,
		such that, for any 
		cube $Q_r(z)$ with radius $r < \rho$ and any
		measurable function $u: Q_r(z) \to \R$, there holds
		\begin{align}\label{comp:eq:lower-bound-balls-abstract}
			\quadForm_{Q_r(z)}^s(u, u)
			&\leq c
			\int_{Q_{\Lambda r}(z)}\int_{Q_{\Lambda r}(z)} \big( u(x) - u(y)\big)^2 \; m(\intd(x, y)).
		\end{align}
	\end{lemma}
	We sketch the strategy here and refer to \cite[Lemma 5.2.9]{HeppDiss} for the full proof.
	The first step is a sort of triangle inequality for the quadratic form $\quadForm^s$ involving cubes in $\Rd$.	
	We introduce the following notation for a cube of radius $R>0$ in $\Rd$.
	\begin{align*}
		Q_j \defeq \bigtimes_{i=1}^d (a_i^{(j)}, b_i^{(j)}) \text{ with }  b_i^{(j)} - a_i^{(j)} = 2R.
	\end{align*}	
	Let $k_1, k_2 \in \{1, \dots d\}$. We say two cubes $Q_1$ and $Q_2$ are \textit{shifted in dimensions $k_1$ and $k_2$} if $(a_i\upperInd{1}, b_i\upperInd{1}) = (a_i\upperInd{2}, b_i\upperInd{2})$ for $i \neq k_1, k_2$. 
	For two shifted cubes $Q_1$ and $Q_2$, we may construct a third cube $Q_3$ of the same radius whose generating intervals are the same as those of $Q_1$ and $Q_2$ in all dimensions $i \neq k_1, k_2$. For the other dimensions, we assume 
	\begin{align*}
		(a_{k_1}\upperInd{3}, b_{k_1}\upperInd{3}) = (a_{k_1}\upperInd{1}, b_{k_1}\upperInd{1}),
		(a_{k_2}\upperInd{3}, b_{k_2}\upperInd{3}) = (a_{k_2}\upperInd{2}, b_{k_2}\upperInd{2}).
	\end{align*}
	Then we can show that there is a constant $c$, depending only on the dimension, such that any measurable function $u$ satisfies 
	\begin{align}\label{comp:eq:goal-for-coercivity-lipschitz-triangle-inequ-for-cubes}
		\int_{Q_1} \int_{Q_2}& U(x, y) |y-x|^{-d-2s} \, \intd y \intd x\\
		&\leq c\; \Big(\int_{Q_1}\int_{Q_3} U(x, y) |y-x|^{-d-2s} \, \intd y \intd x 
		+ \int_{Q_2} \int_{Q_3} U(x, y) |y-x|^{-d-2s} \, \intd y \intd x\Big).\notag
	\end{align}
	The trick to showing this is to construct a cube $Q(|x-y|) \subset Q_3$ for each $x \in Q_1, y \in Q_2$ and multiplying the left-hand side with 
	\begin{align*}
		\frac{1}{|Q(|x-y|)|} \int_{Q(|x - y|)} \intd z = 1,
	\end{align*}
	see \cite[Lemma 5.2.10]{HeppDiss} for details. 
	
	After showing \eqref{comp:eq:goal-for-coercivity-lipschitz-triangle-inequ-for-cubes}, the second step in the proof is to cover any arbitrary cube $Q$ with a set of finitely many disjoint cubes $\mathcal Q$ of some fixed smaller radius, such that $Q \subset \bigcup \mathcal Q \subset \Lambda Q$. Then integrals of the form
	\begin{align*}
		\int_Q \int_Q U(x, y) |y-x|^{-d-2s} \, \intd y \intd x,
	\end{align*}
	are estimated directly using \autoref{comp:ass:comparability-on-small-scales}, while integrals with $Q \neq Q'$ are reduced to
	\begin{align*}
		\int_Q \int_{Q'} U(x, y) |y-x|^{-d-2s} \, \intd x \intd y 
		\leq c \sum_{k=1}^{N(Q, Q')} \int_{Q_k} \int_{Q_k} U(x, y) |y-x|^{-d-2s} \, \intd y \intd x,
	\end{align*}
	using \eqref{comp:eq:goal-for-coercivity-lipschitz-triangle-inequ-for-cubes}. 
	Since $\mathcal Q$ is finite, \eqref{comp:eq:lower-bound-balls-abstract} follows. 
	
	We now prove comparability in Lipschitz domains. The proof works similarly to \autoref{comp:lemma:comparability-in-balls} but uses a more sophisticated covering that is inspired by the proof of \cite[Lemma A.1]{BKS2019}.
	\begin{theorem}\label{comp:thm:lower-bound-for-Lipschitz-domains}
		Let $s \in (0, 1)$ and $\domain \subset \Rd$ be an open, bounded, connected Lipschitz domain. 
		Let $\radAss, \constAss > 0$ and $\lambdaAss \in (0, 1)$. 
		There is a constant $c = c(d, s, \domain, \radAss, \lambdaAss, \constAss) > 0$,
		such that, for any measure $m: \borel(\Rd \times \Rd) \to [0, \infty]$ that satisfies \autoref{comp:ass:comparability-on-small-scales} with radius $\radAss$, factor $\lambdaAss$, and constant $\constAss$, 
		as well as for every measurable function $u: \domain \to \R$, there holds
		\begin{align*}
			\int_\domain \int_\domain \uxy
			\leq c \int_\domain \int_\domain \big(u(x) - u(y)\big)^2 m\big(\intd (x, y) \big). 
		\end{align*}
	\end{theorem}
	
	\begin{proof}
		Let $\domain \subset \Rd$ be a bounded and connected Lipschitz domain
		and $\mathcal W$ a standard Whitney cover of $\domain$ (see \cite[VI.1]{Ste1970}). This means that $\mathcal W$ consists of cubes $Q = Q_r(z)$ with disjoint interiors, such that $2 \sqrt d r \leq \dist(Q, \domain^c) \leq 8 \sqrt d r$. 
		By \cite[Proposition IV.1.2]{Ste1970}, if $\eps \in (0, \frac{1}{4})$,  each point of $\domain$ is contained in at most $3^d$-many enlarged cubes $(1 + \eps)Q$, $Q \in \mathcal W$. The proof of the proposition also shows that $(1 + \eps)Q$ will still be contained in $\domain$. 
		Let $\Lambda' \defeq 1 + \frac{1}{8}$ and $\eta = \eta(d) \defeq \frac{1}{16 \cdot 8 \sqrt d}$.
		Then, if $x \in Q = Q_r(z)$, $\eta \delta_x < \eta 8 \sqrt d r < \frac{1}{16}r$. 
		Setting $\kappa \defeq 1 + \frac{1}{16}$, we see that any $y \in B_{\eta \delta_x}(x)$ also satisfies $y \in Q_{r + \eta \delta_x} (z) \subset  Q_{\kappa r}(z) \subset Q_{\Lambda'r}(z)$.

		We apply inequality (13) in \cite[proof of Theorem 1]{Dyda2006} with $\eta \in (0, 1)$ as chosen above to obtain 
		\begin{align*}
			\int_\domain \int_\domain U(x,y)  |y-x|^{-d-2s} \, \intd y \intd x  
			&\leq c_1 \int_\domain \int_{B_{\eta \delta_x}(x)}U(x,y)  |y-x|^{-d-2s} \, \intd y \intd x  \\
			&\leq c_1 \, \sum_{Q \in \mathcal W} \int_Q \int_{\kappa Q} U(x, y) \, |y-x|^{-d-2s} \, \intd y \intd x.
		\end{align*}
		Now each $Q \in \mathcal W$ is contained in $\domain$, so the radius of each $\kappa Q$ cannot exceed $\kappa \diam(\domain)$. We apply \autoref{comp:lemma:comparability-in-balls} to each $\kappa Q$ with $\Lambda \defeq \frac{\Lambda'}{\kappa} > 1$ and maximal radius $\kappa \diam(\domain)$. The constant $c_2$ we obtain from the lemma will depend on $d, s, \constAss, \lambdaAss, \radAss$, and $\diam(\domain)$, due to the specific choices of $\Lambda$, $\kappa$, and $\rho$. This gives 
		\begin{align*} 
			&\int_\domain \int_\domain U(x,y)  |y-x|^{-d-2s} \, \intd y \intd x 
			\leq c_1 \sum_{Q \in \mathcal W} \int_{\kappa Q} \int_{\kappa Q} U(x, y) \, |y-x|^{-d-2s} \, \intd y \intd x\\
			\leq &\; c_1 c_2 \sum_{Q \in \mathcal W} \int_{\Lambda' Q} \int_{\Lambda' Q}
			U(x, y) \, m \big( \intd (x, y)\big)
			\leq c_1 c_2 3^d \int_\domain \int_\domain U(x, y) \, m \big(\intd x(x, y)\big). 
		\end{align*}
		The last inequality is due to the fact that each point in $\domain$ is only contained in $3^d$-many cubes $\Lambda' Q = \Lambda \kappa Q$. 
	\end{proof}
	
	Note that it is not possible to use \autoref{comp:ass:comparability-on-small-scales} directly in the proof of \autoref{comp:thm:lower-bound-for-Lipschitz-domains} because the enlarging factor $\lambdaAss$ might be too big to allow for the finite overlap property of the Whitney cover. 
	\autoref{thm:coercivity-lip} follows from
	\autoref{comp:thm:lower-bound-for-Lipschitz-domains} together with \autoref{comp:thm:lower-bound-for-singular-measures-small-scales}, which proves \autoref{comp:ass:comparability-on-small-scales} for any singular measure $\frontview(x, \intd y)$ that satisfies the assumptions of \autoref{thm:coercivity-lip}.
	We remark that the argument in \autoref{comp:thm:lower-bound-for-Lipschitz-domains} relies on $\domain$ being bounded, so the coercivity result does not carry over to the entire space $\Rd$. 

\section{Upper bound}\label{sec:upper}
	We will now come to the upper bound in \eqref{eq:comparability}.
	\begin{theorem}[Upper bound for Lipschitz configurations]
		\label{comp:thm:upper-bound}
		Let $s \in (0, 1)$, $\constI > 1$, $\constL, \constW > 0$, $\constExp\in [0, 2s)$, and $\constM \in \N$. 
		Assume $\domain \subset \Rd$ is an open, bounded, and connected Lipschitz domain. 
		There is a constant $c = c (d, s, \domain, \constI, \constL, \constW, \constExp, \constM) > 0$, such that, 
		for any measurable function $u: \domain \to \R$, there holds 
		\begin{align}\label{comp:eq:upper-bound}
			\quadForm^\Gamma_\domain (u, u) \leq c \, \quadForm^s_\domain (u, u), \tag{UB}
		\end{align}
		where $\quadForm^\Gamma_\domain$ has cutoff point $T = \infty$
		and is generated by
		a configuration $\Gamma: \Rd \to (\Sd)^d$ that satisfies 
		\eqref{ass:lip} with Lipschitz constant $\constL$, 
		%\eqref{ass:det} with constant $\constDet \in (0, 1)$, 
		and	\eqref{ass:w} with constants $\constW, \constExp$, and $\constM$, 
		and by a one-dimensional Lebesgue measure $\iota(\intd t)$ that satisfies \eqref{def:iota-comparable} with $\constI$ for $s$. 
	\end{theorem}
	Clearly, \autoref{comp:thm:upper-bound} implies the upper bound in \eqref{eq:comparability} for quadratic forms $\quadForm^\Gamma_\domain$ cut off at any point $T < \infty$. 	
	
	We prove an upper bound for $\quadForm^\Gamma_G(u, u)$, where $G$ is either a small ball in the interior of a domain $\domain$, or a localization box at the boundary of $\domain$, before deducing the inequality \eqref{comp:eq:upper-bound} on the whole domain. 
	If $\domain$ has a Lipschitz boundary, that means there is a localization radius $r_0 > 0$ such that we find, for every point $z \in \partial \domain$, an isometry $T_z: \Rd \to \Rd$ and a Lipschitz function $\gamma_z: \R^{d-1} \to \R$ with Lipschitz constant $\constCov$, such that 
	\begin{align*}
		T_z(\domain) \cap B_{r_0} \big(T_z(z)\big)
		= \{ x \in \Rd \mid x_d > \gamma_z(x')\} \cap B_{r_0}\big(T_z(z)\big).
	\end{align*}
	A \textit{Lipschitz box} with center $T_z(z) \in \partial\domain$, radius $r > 0$, and height $h > 0$ is defined as 
	\begin{align*}
		Q 
		= Q\big(T_z(z), r, h\big) 
		\defeq \{ x \in \Rd \mid x' \in B_r\upperInd{d-1}\big(T_z(z)'\big) \text{ and } 0 < \gamma_z(x') - x_d < h \}. 
	\end{align*}
	Here, $B_r\upperInd{d-1}\big(T_z(z)'\big)$ is a $d-1$-dimensional ball in a $d-1$-dimensional hypersurface of $\Rd$. If we choose $r$ and $h$ small enough with respect to $r_0$, we may always assume $Q \subset \domain$. 
	
	When localizing the quadratic form $\quadForm^\Gamma_\domain (u, u)$ to such a Lipschitz box, the family of measures $\frontview(x, \intd y)$ must also be rotated accordingly. More precisely, if $R_z$ is the rotational part of the isometry $T_z$, then we have to use the measures $\nu^{\Gamma_z}(x, \intd y)$ generated by the configuration $\Gamma_z(x) \defeq \big(R_z(\Gamma_1(T_z\inv(x))), \dots, R_z(\Gamma_d(T_z\inv(x)))\big)$ in the localized integrals over $Q$. 
	However, the results presented here do not depend on any properties of the configuration beside the constants by which they satisfy assumptions \eqref{ass:lip} or \eqref{ass:w}. As the class of configurations satisfying assumptions \eqref{ass:lip} or \eqref{ass:w} with given constants $\constL, \constW, \constExp$, or $\constM$ is closed under isometries, we need not work with the rotated configurations $\Gamma_z$ when dealing with localized integrals over Lipschitz boxes $Q$.
	Thus, in the following, for each localization, we may always assume $\domain$ to be rotated in such a way that $T_z$ is the identity. 
	
	Let us explain a technique to derive upper bounds like 
	\begin{align}\label{comp:eq:localised-upper-bound}
		\int_G \int_G \big(u(x) - u(y)\big)^2 \frontview(x, \intd y) \intd x
		\leq c \int_\domain \int_\domain \big(u(x) - u(y)\big)^2 |y-x|^{-d-2s} \intd y \intd x, 
	\end{align}
	where $G$ is either a ball or a Lipschitz box at the boundary, and $\frontview(x, \intd y)$ is given in \autoref{def:frontview-concrete}. 
	The main idea is to add an "empty" integral over a ball $B_\rho\big(p(x, t)\big)$, where the center $p(x, t)$ and the radius $\rho$ are chosen in such a way that $\big(u(x) - u(x + t \Gamma_k(x))\big)^2 |t|^{-1-2s}$ can be estimated from above by $\big(u(x) - u(y)\big)^2 |y-x|^{-d-2s} + \big(u(y) - u(x + t \Gamma_k(x))\big)^2 |y - (x + t \Gamma(x))|^{-d-2s}$, for any $y \in B_\rho\big(p(x, t)\big).$ 
	
	First, let us observe that it is always enough to prove \eqref{comp:eq:localised-upper-bound} for $\frontview(x, \intd y)$ being generated by one subconfiguration that satisfies the given assumptions. So, from now on, let $\Gamma: \Rd \to \Sd$ be a subconfiguration. 
	\begin{remark}\label{rem:view-map-bi-lipschitz}
		For any $\constL > 0$, there is a radius $R = R (d, \constL) > 0$, such that, for each $|t| < R$, the path map $\viewMap_t$ generated by any subconfiguration $\Gamma: \Rd \to \Sd$ that satisfies \eqref{ass:lip} with $\constL> 0$ is injective on any ball of radius at most $R$. This is due to \autoref{measures:lemma:properties-of-det-view-map} (b).
	\end{remark}
	Now let $r_0$ be some given maximal radius. In practice, we will always be working in some bounded domain, so assuming a maximal radius is no great restriction. 
	Assume also that we are given some enlargement constant $\lambda > 1$. 
	Choose $r < \frac{R}{2 \lambda} \wedge r_0$, where $R$ is the radius from \autoref{rem:view-map-bi-lipschitz}. 
	Let $G$ be either a ball $B_r(z)$ or a Lipschitz box $Q \subset B_r(z) \cap \domain$, such that $\tilde G$ is either $B_{\lambda r}(z)$ or $B_{\lambda r}(z) \cap \domain$. 
	The choice of $r$ ensures $|t| < 2 \lambda r_0 \wedge R$ if $x, x+ t \Gamma(x) \in \tilde G$. 
	In particular $\viewMap_t$ will be injective on $\tilde G$ in this case. 
	We assume
	\begin{align}\label{comp:eq:fatness-techinque-1}
		\exists \, c_1, c_2 > 0, \rho: [0, \infty) \to [0, \infty) \text{ measurable, s.t. } \forall \; \tau \in [0, R] 
		\; : \; 
		c_1 \tau \leq \rho(\tau) \leq c_2 \tau.\tag{U1}
	\end{align}
	Assume furthermore that the center $p(x, t)$ is chosen such that 
	\begin{align}\label{comp:eq:fatness-techinque-2}
		\exists \; c_3 > 0 \; 
		\forall \; x \in G, t \in (-2r, 2r) \; : \; 
		\big|p(x, t) -x\big|, \big|p(x, t) - \big(x + t \Gamma(x)\big)\big|
		\leq c_3 |t|. \tag{U2}
	\end{align}
	For the entire ball, assume 
	\begin{align}\label{comp:eq:fatness-techinque-3}
		\forall \; x \in G, t \in (-2r, 2r) \; : \; B_\rho\big(p(x, t)\big) \subset \tilde G.\tag{U3}
	\end{align}
	Then there are constants $c' =  c'(d, c_1\inv) > 0$, such that
	\begin{align}\label{comp:eq:fatness-techniwue-integral}
		1 \leq c' |t|^{-d} \int_{B_\rho(p(x, t))} \intd y 
		\text{ where }
		c' \defeq \frac{1}{c_1^d |B_1|}.\tag{V1}
	\end{align}
	and $c''(c_2, c_3) > 0$ such that, for any $y \in B_\rho\big(p(x, t)\big)$, there holds 
	\begin{align}\label{comp:eq:fatness-technique-norm}
		|(t, y-z)| \leq c'' |t| \text{ where } 
		(t, y-z) \defeq 
		\begin{pmatrix}
			t\\
			y_1 - z_1\\
			\vdots\\
			y_d - z_d
		\end{pmatrix}, 
		\text{ and } z \in \{x, x+t\Gamma(x)\}.\tag{V2}
	\end{align}
	Inequality \eqref{comp:eq:fatness-techniwue-integral} follows directly from \eqref{comp:eq:fatness-techinque-1}, and \eqref{comp:eq:fatness-technique-norm} from \eqref{comp:eq:fatness-techinque-2}. 
	
	Let us prove \eqref{comp:eq:localised-upper-bound} under assumption \eqref{comp:eq:fatness-techinque-1} - \eqref{comp:eq:fatness-techinque-3}.
	We first use \eqref{comp:eq:fatness-techniwue-integral} to add the additional integral over $B_\rho\big(p(x, t)\big)$. Note that the function $(x, t) \mapsto \int_{B_{\rho(|t|)}(p(x, t))}\intd y = \rho(|t|)^d |B_1|$ is measurable by \eqref{comp:eq:fatness-techinque-1}. After applying \eqref{def:iota-comparable}, the estimate \eqref{comp:eq:fatness-techniwue-integral} changes the exponent of $|t|$ from $-1-2s$ to $-d-1-2s$. 
	Then we insert an additional point $y \in B_\rho\big(p(x,t)\big)$ and apply \eqref{comp:eq:fatness-technique-norm} to $|t|^{-d-1-2s}$, 
	\begin{align*}
		&\int_{G} \int_\R 
		U\big(x, x + t \Gamma(x)\big)
		\indic_G\big(x + t \Gamma(x)\big)
		|\det D \viewMap_t(x)| 
		\; \iota(\intd t) \intd x\\
		\leq &\; c'\, \constI \int_{G} \int_\R \int_{B_\rho(p(x, t))}
		U\big(x, x+t\Gamma(x)\big)
		\indic_G\big(x + t \Gamma(x)\big)
		%&\phantom{\; c' \int_{G} \int_\R \int_{B_\rho(p(x, t))}}
		| \det D \viewMap_t(x)|
		\,|t|^{-d-1-2s}
		\; \intd y \intd t \intd x\\
		\leq &\; c(c', c'', \constI, d, \constL,  r_0, \lambda) \int_{G} \int_\R \int_{B_\rho(p(x, t))}
		U(x, y)
		|(t, y-x)|^{-d-1-2s}
		\; \intd y \intd t \intd x\\
		& + \; c( c', c'', \constI) \int_{G} \int_\R \int_{B_\rho(p(x, t))}
		U\big(y, x + t \Gamma(x)\big)
		\indic_G\big(x + t \Gamma(x)\big)\\
		&\phantom{mmmmmmmmmmmmmmmm}
		| \det D \viewMap_t(x)|
		\,|(t, y-(x + t \Gamma(x)))|^{-d-1-2s}
		\; \intd y \intd t \intd x.
	\end{align*}
	We have also used \autoref{measures:lemma:properties-of-det-view-map} (a) to estimate $|\det D \viewMap_t(x)|$ by a constant $c (d, \constL, r_0, \lambda)$, and estimated $\indic_G$ by 1 in the first summand. 
	%	Observe 
	%	\begin{align}\label{comp:eq:barteks-trick}
		%		\int_\R \big|(t, y-x)\big|^{-d-1-2s} \intd t 
		%		&= |y-x|^{-d-1-2s} \int_\R \Big( \frac{t^2}{|y-x|^2} + 1\Big)^{-\frac{d + 1 + 2s}{2}} \intd t\notag\\
		%		&\leq c(d) |y-x|^{-d-2s},
		%	\end{align}
	%	by a change of variables $\tau(t) \defeq t/|y-x|$ (c.f. \cite[Proof of Theorem 6.2]{DyKa20}). 
	Integrating out $t$ as before, and applying \eqref{comp:eq:fatness-techinque-3} to the first summand gives
	$\int_{G} \int_{\tilde G} U(x, y) |y-x|^{-d-2s} \; \intd y \intd x$.
	For the second summand, we first use Tonelli and \eqref{comp:eq:fatness-techinque-3}, before we use the are formula with $z(x) \defeq V_t(x)$, keeping in mind that $\viewMap_t$ is injective on $G$ due to the choice of $r$:
	\begin{align*}
		&\int_{\tilde G} \int_{-2R}^{2 R} \int_\Rd
		U\big(y, \viewMap_t(x)\big)
		\indic_G(\viewMap_t(x))
		| \det D \viewMap_t(x)|
		|(t, y-\viewMap_t(x))|^{-d-1-2s}
		\; \intd x \intd t \intd y\\
		= & \; \int_{\tilde G} \int_{-2R}^{2R} \int_G
		U(y, z)
		|(t, y-z)|^{-d-1-2s}
		\; \intd z \intd t \intd y 
	\end{align*}
	Then integrating out $t$ leads us to
	$\int_{\tilde G} \int_G U(x, y) |y-x|^{-d-2s} \, \intd y \intd x$ as before.
	Overall, we obtain 
	\begin{align*}
		&\int_G \int_G U(x, y) \, \frontview(x, \intd y) \intd x
		\leq \; c\; \int_{\tilde G} \int_{\tilde G} U(x, y) \,|y-x|^{-d-2s} \,\intd y \intd x,
	\end{align*}
	where the constant depends on $d, s, \domain, \constL$ as well as $c_1, c_2,$ and $c_3$. 
	%The bound \eqref{comp:eq:localised-upper-bound} is true for families of measures $\frontview(x, \intd y)$ with cutoff point $T = \infty$. This immediately implies \eqref{comp:eq:localised-upper-bound} for families of measures with cutoff point $T < \infty$. 
	
	Now we can prove \eqref{comp:eq:localised-upper-bound} for balls and Lipschitz boxes. By the argument just presented, it will suffice to verify \eqref{comp:eq:fatness-techinque-1} - \eqref{comp:eq:fatness-techinque-3} in both cases. Since both proofs are not hard to perform, we will merely state the respective results. 
	For balls, with have (see also \cite[Lemma 5.3.3]{HeppDiss}):  
	\begin{lemma}\label{comp:lemma:upper-bound-small-scales-balls}
		Let $s \in (0, 1)$, $\constL, r_0> 0$, and $\lambda, \constI > 1$. 
		There are a radius $R_0 = R_0 (d, \constL) \in (0, r_0)$ and 
		a constant $c = c(d, s, \constI, \constL, r_0, \lambda) > 0$, such that,
		for any ball $B= B_R(z)$ with $R \in (0, R_0)$, and 
		any measurable function $u : \lambda B \to \R$, there holds 
		\begin{align*}
			\quadForm^\Gamma_{B} (u, u) 
			\, \leq \; c \, \quadForm^s_{\lambda B}(u, u),
		\end{align*}
		where $\quadForm^\Gamma_B$ is generated by 
		any configuration $\Gamma: \Rd \to (\Sd)^d$ 
		satisfying \eqref{ass:lip} with Lipschitz constant $\constL$, 
		and any one-dimensional Lebesgue measure $\iota(\intd t)$ satisfying \eqref{def:iota-comparable} with $\constI$ for $s$. 
	\end{lemma}
	The proof shows that $c \to \infty$ as $\lambda \to 1$, which makes it impossible to obtain, by this technique, a result where the ball does not grow at least slightly by going to the limit. 
	%	\begin{proof} 
	%		Let $R_0 \defeq R_1/2 \wedge r_0$, where $R_1$ is the radius from \autoref{rem:view-map-bi-lipschitz}. 
	%		Let $\Gamma: \Rd \to \Sd$ be a subconfiguration with Lipschitz constant $\constL$ and assume $x, x + t \Gamma(x) \in B_r(z)$ where $r < R_0$. This implies $|t| \in [0, 2r)$. Then $\rho(\tau) \defeq \tau \wedge (1-\lambda)r$ satisfies \eqref{comp:eq:fatness-techinque-1} with $c_1 = (1-\lambda)/2$ and $c_2 = 1$. 
	%		With $p(x, t) \defeq x$, \eqref{comp:eq:fatness-techinque-2} is satisfied as well for $c_3 = 5$. 
	%		As $B_\rho(x) \subset B_{\lambda r}(z)$ for any $x \in B_r(z)$, we obtain \eqref{comp:eq:fatness-techinque-3}. Thus, the conclusion follows with the technique explained above. 
	%	\end{proof}
	%Now any minimal distance $\delta > 0$ to the boundary of $\domain$ allows us to choose $R_0 < \delta / \lambda$, such that, if $\dist(z, \domain^c) > \delta$, there holds $B_r(z) \subset B_{\lambda r}(z) \subset \domain$ for any $r \in (0, R_0)$. This implies \eqref{comp:eq:localised-upper-bound} for $m \big(\intd (x, y)\big) = |y - x|^{-d-2s} \intd y \intd x$. 
	
	The lemma for Lipschitz boxes reads (c.f. \cite[Lemma 5.2.6]{HeppDiss}): 
%	\eqref{comp:eq:localised-upper-bound} for Lipschitz boxes. 
%	In a Lipschitz box with Lipschitz function $\gamma_z$ and Lipschitz constant $\constCov$, the \textit{$\constCov$-cone at $x$} is given by
%	\begin{align*}%\label{comp:def:lipschitz-cone}
%		V_\constCov[x] \defeq
%		\{ h \in \Rd \mid h_d - x_d > \constCov|h'-x'| \}.
%	\end{align*}
%	This cone lies above the graph $\Gamma_{\gamma_z}$, i.e. $h \in V_\constCov[x]$ and $h' \in B_r\upperInd{d-1}(z')$ together imply $\gamma(h') < h_d$.
%	Consequently, the distance to the boundary of any point in such a cone may be estimated from below by its distance to the boundary of the cone. More precisely, if $p \in V_\constCov[x]$, then 
%	\begin{align*}
%		\dist(p, V_\constCov[x]^c)
%		= \frac{(p_d - x_d) - \constCov|p' - x'|}{\sqrt{1+\constCov^2}}.
%	\end{align*}
%	The proof of this latter fact can also be found in Step 2 of the proof of \cite[Proposition 5]{Dyda2006}. For more details, see \cite[Lemmas 5.3.4. and 5.3.5]{HeppDiss}.
%	
%	Using these observations, we show the upper bound for Lipschitz boxes. 
	\begin{lemma}\label{comp:lemma:upper-bound-lipschitz-box}
		Let $\domain \subset \Rd$ be an open and bounded Lipschitz domain with localization radius $r_0 > 0$ and Lipschitz constant $\constCov > 0$. 
		Let $s \in (0, 1)$, $\constI>1$, and $\constL > 0$. 
		There are a 
		radius $R_0 (d, \constL, \constCov, r_0) \in (0, r_0)$
		and a constant $c = c(d, s, \diamDom, \constI, \constL, \constCov) > 0$, 
		such that, for any $r \in (0, R_0)$, any Lipschitz box $Q$ of $\domain$ at $z \in \partial\domain$ with $Q \subset B_{r}(z)$, and any measurable function $u: \domain \to \R$, there holds 
		\begin{align*}
			\quadForm^\Gamma_Q (u, u) 
			\, \leq \; c \, \quadForm^s_\domain(u, u),
		\end{align*} 
		where $\quadForm^\Gamma_Q$ is generated by 
		any configuration $\Gamma: \Rd \to (\Sd)^d$ 
		satisfying \eqref{ass:lip} with Lipschitz constant $\constL$, 
		and any one-dimensional Lebesgue measure $\iota(\intd t)$ satisfying \eqref{def:iota-comparable} with $\constI$ for $s$. 
	\end{lemma}

	Let us now tackle \autoref{comp:thm:upper-bound}.
	The proof is heavily inspired by the proof of Theorem 1.c in \cite{Dyda2006}. The main idea is to bound the integral over $G \times G'$, where $G$ and $G'$ are two distinct elements in the cover of $\domain$, by a sum of integrals over sets $G_i \times G_i$, which is expressed in the following lemma: 	
	\begin{lemma}\label{comp:lemma:upper-bound-lipschitz-domain-laplace-localized}
		Let $s \in (0, 1)$, $\constI > 1$, $\constL, \constW > 0$, $\constExp\in [0, 2s)$, $\constM \in \N$, and  $\domain \subset \Rd$ be an open, bounded, and connected Lipschitz domain. 
		There is a constant $c = c (d, s, \domain, \constI, \constL, \constW, \constExp, \constM) > 0$, such that, 
		for any finite cover $\mathcal G$ of $\domain$ that consists of Lipschitz boxes at the boundary and balls in the interior of $\domain$ and satisfies 
		\begin{align}\label{comp:eq:cover-overlapping-property-boundary}
			\exists \; \eps > 0 \; \forall \; z \in \partial \domain \;
			\exists \;Q \in \mathcal G \; :\;
			B_\eps(z) \cap \domain \subset Q,
		\end{align}
		as well as 	
		for any measurable function $u: \domain \to \R$, there holds 
		\begin{align*}
			\quadForm^\Gamma_\domain (u, u) 
			\leq c \sum_{G \in \mathcal G} \quadForm^s_G (u, u), 
		\end{align*}
		where $\quadForm^\Gamma_\domain$ has cutoff point $T = \infty$
		and is generated by
		a configuration $\Gamma: \Rd \to (\Sd)^d$ that satisfies 
		\eqref{ass:lip} with Lipschitz constant $\constL$, 
		%\eqref{ass:det} with constant $\constDet \in (0, 1)$, 
		and	\eqref{ass:w} with constants $\constW, \constExp$, and $\constM$, 
		and by a one-dimensional Lebesgue measure $\iota(\intd t)$ that satisfies \eqref{def:iota-comparable} with $\constI$ for $s$. 
	\end{lemma}
	A proof of this lemma can be found in \cite[Lemma 5.3.7]{HeppDiss}. 
	\begin{proof}[Proof of \autoref{comp:thm:upper-bound}]
		Let $\domain \subset \Rd$ be an open, bounded, and connected Lipschitz domain. Let $R_Q$ be the maximal radius from \autoref{comp:lemma:upper-bound-lipschitz-box} for Lipschitz boxes. Since $\domain$ is bounded, we may cover $\partial \domain$ with a finite set of Lipschitz boxes $\mathcal Q$ such that \eqref{comp:eq:cover-overlapping-property-boundary} holds. 
		
		Now we cover the interior with balls. Set $\domain_{\delta} \defeq \{ x \in \domain \mid \dist(x, \partial \domain) > \delta\}$ for $\delta > 0$. 
		There will be some small $\delta_0 > 0$, such that $\domain \setminus \domain_{\delta_0} 
		\subset \bigcup_{Q \in \mathcal Q} Q,$
		i.e. the Lipschitz boxes cover a little more than just the boundary of $\domain$. 
		Applying \autoref{comp:lemma:upper-bound-small-scales-balls} with $\lambda = 2$ and radius bounded by $\diamDom$, we find a radius $R_B < \delta_0/2$, such that  \eqref{comp:eq:localised-upper-bound} holds for any ball $B_{R_B}(x)$ with $x \in \domain_{\delta_0}$. Cover $\domain_{\delta_0}$ with balls of radius $R_B$. This cover, too, will be finite, due to the boundedness of $\domain$. 
		
		Let $\mathcal G$ be the set that contains all Lipschitz boxes and balls in the cover of $\domain$.
		$\mathcal G$ was chosen such that either \autoref{comp:lemma:upper-bound-small-scales-balls} or \autoref{comp:lemma:upper-bound-lipschitz-box} holds for each of its members. 
		Together with \autoref{comp:lemma:upper-bound-lipschitz-domain-laplace-localized}, this implies
		\begin{align*} 
			\quadForm_\domain^\Gamma(u, u) 
			\, \leq \, c \, \sum_{G \in \mathcal G} \quadForm_G^s(u, u) 
			\, \leq \, c \, \quadForm^s_\domain(u, u).
		\end{align*}
	\end{proof}
	
	\begin{small}
	\bibliographystyle{plain}%alpha
	\bibliography{lit}
	\end{small}

	\bigskip
	\bigskip
	
	\noindent
	\textsc{Fakultät für Mathematik, Universität Bielefeld \\Postfach 100131, D-33501 Bielefeld}
	
	\noindent
	Kai-Uwe Bux: 
	\texttt{bux@math.uni-bielefeld.de}
	
	\noindent
	Solveig Hepp: 
	\texttt{shepp@math.uni-bielefeld.de}

\end{document}